\documentclass[11pt,english]{amsart}

\def\margin_comment#1{\marginpar{\sffamily{\tiny #1\par}\normalfont}}
\usepackage[table]{xcolor}
\usepackage{amssymb,euscript}
\usepackage{amsfonts}
\usepackage{geometry}
\usepackage{fancyhdr}
\usepackage{graphicx}
\usepackage{setspace}
\usepackage{color}
\usepackage{amsmath,amssymb,amsthm}
\usepackage{epsfig}
\usepackage{babel}
\usepackage{tikz}
 \usepackage{amsmath}
\usepackage{tikz-cd}
\usetikzlibrary{matrix, calc, arrows}
\usepackage{tikz}
\usetikzlibrary{arrows,chains,matrix,positioning,scopes}
\usepackage{lipsum}
\usepackage{mathtools}
\usepackage{mathrsfs}
\usepackage{graphicx}
\usepackage{graphics}
\usepackage{setspace}
	\usepackage{color}

\usetikzlibrary{trees}
\usetikzlibrary{arrows,shapes,automata,backgrounds,petri,positioning,scopes,matrix, calc, plothandlers,plotmarks,patterns}
\usepackage{enumerate}
\usepackage{float}
\usepackage[utf8]{inputenc}
\usepackage[T1]{fontenc}
\usepackage[document]{ragged2e}
 \usepackage{verbatim}
\usepackage{tikz-3dplot}
\usepackage{pgfplots}
\usepackage{tkz-graph}
\GraphInit[vstyle = Normal]
\tikzset{
  LabelStyle/.style = {minimum width = 2em, 
                        text = red, font = \bfseries },
  VertexStyle/.append style = { inner sep=2pt,
                                font = \Large\bfseries, fill},
  EdgeStyle/.append style = {->, bend left} }

\makeatletter
\tikzset{join/.code=\tikzset{after node path={%
\ifx\tikzchainprevious\pgfutil@empty\else(\tikzchainprevious)%
edge[every join]#1(\tikzchaincurrent)\fi}}}
\makeatother
\tikzset{>=stealth',every on chain/.append style={join},
         every join/.style={->}}
\tikzstyle{labeled}=[execute at begin node=$\scriptstyle,
   execute at end node=$]
\newtheorem{thm}{Theorem}[section]
\numberwithin{equation}{section} 
\numberwithin{figure}{section} 
\theoremstyle{plain}
\newtheorem*{thm*}{Theorem}
\theoremstyle{definition}
\theoremstyle{plain}
\newtheorem{thm_A}{Theorem}
\newtheorem*{defn*}{Definition}

\theoremstyle{plain}

\theoremstyle{plain} 

\theoremstyle{plain}
\newtheorem{cor_A}{Corollary}
\newtheorem{prop}[thm]{Proposition} 
\theoremstyle{remark}
\newtheorem{ex}[thm]{Example}
\theoremstyle{remark}
\newtheorem{rem}[thm]{Remark}
\theoremstyle{plain}

\theoremstyle{plain}

\theoremstyle{plain}
\newtheorem{lem}[thm]{Lemma} 
\theoremstyle{definition}
\newtheorem{defn}[thm]{Definition}
\newtheorem*{acknowledgment}{Acknowledgment}
\newtheorem*{acknowledgment*}{Addentum}

\theoremstyle{plain}
\newtheorem*{ex*}{Example}
\theoremstyle{plain}

\begin{document}
\title[The Yang-Baxter equation and cryptography ]{The Yang-Baxter equation  and cryptography }
\author{Fabienne Chouraqui}
\begin{abstract}
	We find a method to construct iteratively  from  a non-degenerate involutive  set-theoretic solution of 	the  Yang-Baxter equation an infinite family of very large non-degenerate involutive  set-theoretic solutions. In case the initial solution is irretractable, all the induced  solutions are also irretractable. In case the initial solution is indecomposable, we give a criterion to decide whether all the induced  solutions are also indecomposable.  Besides the interest in the construction of  large (indecomposable) solutions of the Yang-Baxter equation, this construction may have some applications in cryptography. Indeed,  we suggest a public key encryption method and a signature method based on our construction,  and examine  their strengths and weaknesses.
\end{abstract}
\maketitle

\section*{Introduction}
The Yang-Baxter equation is an equation in mathematical physics and it lies in the  foundation of  the theory of quantum groups. One of the fundamental problems is to find all the solutions of this equation. In \cite{drinf}, Drinfeld suggested the study of a particular class of solutions,  derived from the so-called set-theoretic solutions.  A  set-theoretic solution of the Yang-Baxter equation is a pair $(X,r)$, where $X$ is a set and 
\[r: X \times X \rightarrow X \times X\,,\;\;\; r(x,y)=(\sigma_x(y),\gamma_y(x))\]
is a bijective map satisfying $r^{12}r^{23}r^{12}=r^{23}r^{12}r^{23}$, where $r^{12}=r \times Id_X$ and  $r^{23}=Id_X\times r$. A set-theoretic solution  $(X,r)$ is said to be non-degenerate if, for every $x \in X$, the maps $\sigma_x,\gamma_x$ are bijections of $X$ and it  is said to be involutive if $r^2=Id_{X\times X}$.   Non-degenerate and involutive set-theoretic solutions give rise to solutions of the  Yang Baxter.  Indeed, by defining  $V$ to be  the  real vector space spanned by $X$,  and $c:V \otimes V \rightarrow V \otimes V$  to be the  linear operator induced by  $r $, then   $c$ is  a linear operator satisfying the equality  $c^{12}c^{23}c^{12}=c^{23}c^{12}c^{23}$ in $V \otimes V \otimes V$, that is $c$  is a solution of the  Yang-Baxter equation. By composing $r$ with $\tau$, 
where $\tau$ is the flip map $\tau(x,y)=(y,x)$,  the induced linear operator, $R$,  is  a linear operator satisfying the equality  $R^{12}R^{13}R^{23}=R^{23}R^{13}R^{12}$ in $V \otimes V \otimes V$, that is $R$  is a solution of the quantum Yang-Baxter equation. 


Another source of solutions of the  Yang-Baxter equation  is  derived from quasi-triangular bialgebras (or braided bialgebras). In \cite{drinf-quasi1,drinf-quasi2}, Drinfeld introduced the concept of quasi-triangular bialgebras. These are bialgebras with a  universal R-matrix inducing a solution of  the   Yang-Baxter equation on any of their modules. There is a dual notion of cobraided bialgebras, these are bialgebras inducing  $R$-matrices on their comodules.   Both constructions provide a systematic method to produce solutions of   the  Yang-Baxter equation. Conversely, given a  solution of  the  Yang-Baxter equation, 
there is a method, due to Faddeev, Reshetikhin and  Takhtadjian (FRT construction for short),  to construct a cobraided bialgebra  \cite{frt}, \cite{kassel}. Faddeev, Reshetikhin and Takhtadjian show that if  $c \in \operatorname{End}(V \otimes V)$, where  $V$  is a finite-dimensional vector space, is a solution of the  Yang-Baxter equation, then there exists  a cobraided bialgebra $A(c)$ coacting on $V$ \cite{frt}. 
The algebra  $A(c)$ is the quotient of  the free algebra  by a two-sided  ideal $I(c)$. We  consider the special case,  where  $c$ is a solution of the YBE induced by  a non-degenerate and involutive solution $(X,r)$, with $\mid X\mid=n$, and compute the set  of generating elements of the two-sided ideal $I(c)$ in  $A(c)$.  By considering  only this set  of generating elements  and forgetting about the  algebraic structures  $A(c)$ and  $I(c)$,  we find a method to  define  a new non-degenerate and involutive solution of size $n^2$, and then iteratively construct a  non-degenerate and involutive solution, $(\mathcal{X},  \tilde{r})$,  of size $n^{2^k}$,  for every natural number $k \geq 1$.   That is, we prove the following Theorem,   with all the precise definitions  given later:
 \begin{thm_A}\label{thmA-construction}
 	Let $(X,r)$ be a  non-degenerate and involutive set-theoretic solution of the   Yang-Baxter equation, with $\mid X\mid=n$.  	Then, for every natural number $k \geq 1$,    there exists a  non-degenerate and involutive set-theoretic solution, $(\mathcal{X},  \tilde{r})$,  of size $n^{2^k}$ induced from $(X,r)$. Furthermore, the following hold:
 	\begin{enumerate}[(i)]
 		\item If  $(X,r)$  is irretractable or a multipermutation solution of level $\ell$, then 
 		$(\mathcal{X},  \tilde{r})$ is  irretractable or  a multipermutation solution of level $\ell$, respectively.
 			\item 	If  $(X,r)$  is decomposable, then 
 		$(\mathcal{X},  \tilde{r})$ is also decomposable.
 		\item 	If  $(X,r)$  is of class $m$, then 
 		$(\mathcal{X},  \tilde{r})$ is also  of class $m$.
 		 	\item 	If  $(X,r)$  is 	indecomposable and satisfies  a certain  condition $( \mathfrak{C})$, then 	$(\mathcal{X},  \tilde{r})$  is indecomposable and satisfies $( \mathfrak{C})$.
 	\end{enumerate}
 \end{thm_A}
It holds also that if  $(X,r)$ and $(X',r')$ are different non-degenerate and involutive set-theoretic solutions of the   Yang-Baxter equation, with $\mid X\mid=\mid X'\mid=n$, then their induced solutions $(\mathcal{X},  \tilde{r})$  and $(\mathcal{X}',  \tilde{r}')$  are different, and if  $(X,r)$ and $(X',r')$ are isomorphic, 	then $(\mathcal{X},  \tilde{r})$  and $(\mathcal{X}',  \tilde{r}')$ are also isomorphic.
\begin{cor_A}\label{cor1}
	Let $(X,r)$ be a  non-degenerate and involutive set-theoretic solution of the   Yang-Baxter equation, with $\mid X\mid=n$.  If  $(X,r)$  is 	irretractable, indecomposable and satisfies  a certain  condition $( \mathfrak{C})$,  
	then for every natural number $k \geq 1$,    there exists a  non-degenerate and involutive set-theoretic solution of size $n^{2^k}$, that is also  irretractable and indecomposable.
\end{cor_A}

Based on our method of construction of huge  solutions we suggest the idea that there may be some  possible applications  to cryptography. We suggest a public key encryption method and a signature method based on our construction, we evaluate their complexity and examine their strengths and weaknesses.

 The paper is organized as follows. In Section $1$, we give some preliminaries on set-theoretic solutions of the Yang-Baxter equation, and the FRT construction. In Section $2$, we give some preliminaries on  cryptography. In Section $3$, we  present the method of construction of large  set-theoretic solutions and we prove that the properties  non-degenerate and involutive are inherited by the induced solutions. In the particular case that the initial solution is indecomposable, we give a criterion to decide whether all the induced  solutions are also indecomposable  and  prove Theorem $1$. In Section $4$, we present the possible applications  of our method of construction of huge  solutions to cryptography. We suggest a public key encryption method and a signature method based on our construction, we evaluate their complexity and examine their strengths and weaknesses. In the appendix, we suggest a tentative key exchange method based on our construction.
 \begin{acknowledgment}
 	I am very grateful to Leandro Vendramin   for his  great help in learning the gap package "YangBaxter" and for the construction of Table \ref{table_nb-solutions}.
 \end{acknowledgment}
The data that support some of the findings of this study (like the number and the enumeration of solutions of size less than 10) are openly available in  the gap package "YangBaxter".
\section{Preliminaries on set-theoretic solutions of the   Yang-Baxter equation  (YBE) }
\subsection{Definition and properties of set-theoretic solutions of the YBE}
\label{subsec_qybe_Backgd} There is a vast literature  on the topic and we  refer to  \cite{catino3}, \cite{cedo}, \cite{jespers-adv}, \cite{brace}, \cite{etingof}, \cite{gateva_van,
	gateva_new},  \cite{g-vendra1}, \cite{jespers_book},  \cite{lebed1}, \cite{adolfo}, \cite{rump_braces}, \cite{smot} and many more.
Let $X$ be a non-empty set. Let $r: X \times X \rightarrow X \times X$  be a map and write $r(x,y)=(\sigma_{x}(y),\gamma_{y}(x))$,  where $\sigma_x, \gamma_x:X\to X$ are functions  for all  $x,y \in X$.   The pair $(X,r)$ is  \emph{braided} if $r^{12}r^{23}r^{12}=r^{23}r^{12}r^{23}$, where the map $r^{ii+1}$ means $r$ acting on the $i$-th and $(i+1)$-th components of $X^3$.  In this case, we  call  $(X,r)$  \emph{a set-theoretic solution of the  Yang-Baxter equation}, and whenever $X$ is finite, we  call  $(X,r)$  \emph{a finite set-theoretic solution of the  Yang-Baxter equation}.  The pair $(X,r)$ is \emph{non-degenerate} if for every  $x\in X$,  $\sigma_{x}$ and $\gamma_{x}$  are bijective and it   is  \emph{involutive} if $r\circ r = Id_{X^2}$. If $(X,r)$ is a non-degenerate involutive set-theoretic solution, then $r(x,y)$ can be described as  $r(x,y)=(\sigma_{x}(y),\gamma_{y}(x))=(\sigma_{x}(y),\,\sigma^{-1}_{\sigma_{x}(y)}(x))$.  A set-theoretic solution  $(X,r)$ is \emph{square-free},  if  for every $x \in X$, $r(x,x)=(x,x)$. A set-theoretic  solution $(X,r)$ is \emph{trivial} if $\sigma_{x}=\gamma_{x}=Id_X$, for every  $x \in X$. 
 \begin{lem}\cite{etingof} \label{lem-formules-invol+braided}
 	\begin{enumerate}[(i)]
 	\item  $(X,r)$ is involutive  if  and only if  for every  $x,y \in X$:
 	\begin{gather}
 \sigma_{\sigma_x(y)}\gamma_{y}(x)=x \label{eqn-inv1}\\
 \gamma_{\gamma_y(x)}\sigma_x(y)=y \label{eqn-inv2}
 	\end{gather} 
\item   $(X,r)$ is  braided if  and only if, 	for  every $x,y,z \in X$, the following  holds:
 \begin{gather}
\sigma_x\sigma_y=\sigma_{\sigma_x(y)}\sigma_{\gamma_y(x)} \label{eqn-braided-sigma}\\  \gamma_y\gamma_x=\gamma_{\gamma_y(x)}\gamma_{\sigma_x(y)} \label{eqn-braided-gamma}\\\gamma_{\sigma_{\gamma_y(x)}(z)}(\sigma_x(y))=\sigma_{\gamma_{\sigma_y(z)}(x)}(\gamma_z(y)) \label{eqn-braided-old}
\end{gather}
 \end{enumerate}
\end{lem}

\begin{defn}
	Let  $(X,r)$   be a set-theoretic solution of the YBE. The \emph{structure group} of $(X,r)$ is  defined by $G(X,r)=\operatorname{Gp} \langle X\mid\ xy =\sigma_x(y)\gamma_y(x)\ ;\ x,y\in X \rangle$.  
\end{defn}
The  structure group of the trivial solution is $\mathbb{Z}^{X }$.  Two  set-theoretic solutions $(X,r)$ and $(X',r')$ are \emph{isomorphic} if there is a bijection $\mu:X \rightarrow X'$ such that $(\mu \times \mu) \circ r=r'\circ (\mu \times \mu)$  \cite{etingof}. If  $(X,r)$ and $(X',r')$ are isomorphic, then $G(X,r) \simeq G(X',r')$, with  $G(X,r)$ and $G(X',r')$ their respective structure groups. 
\begin{defn}\cite{cedo}
		Let  $(X,r)$   be a set-theoretic solution of the YBE.
	The subgroup of $\operatorname{Sym}_X$  generated by $\{\sigma_x\mid x\in X\}$ is denoted by $\mathcal{G}(X,r)$ and is called \emph{a IYB group}. 
\end{defn}

\begin{defn}
Let  $(X,r)$ be a non-degenerate and involutive set-theoretic solution of the YBE. 
\begin{enumerate}[(i)]
	\item  A set $Y \subset X$ is  \emph{invariant} if  $r(Y \times Y)\subset Y \times Y$.
	\item An invariant  subset $Y \subset X$  is  \emph{non-degenerate} if  $(Y,r\mid_{Y^2})$ is non-degenerate involutive set-theoretic solution of the YBE. 
	\item  $(X,r)$  is \emph{decomposable} if it is a union of two non-empty disjoint non-degenerate invariant subsets. Otherwise, it is called \emph{indecomposable}.
\end{enumerate}
\end{defn}
A very simple class of  non-degenerate involutive set-theoretic solutions of the YBE is the class of  \emph{permutation solutions}. These solutions have the form $r(x,y)=(\sigma(y),\sigma^{-1}(x))$,  where the bijections  $\sigma_x: X\to X$  are  all equal and equal to $\sigma$, the   bijections  $\gamma_x: X\to X$  are  all equal and equal to $\sigma^{-1}$. If $\sigma$ is a cyclic permutation, $(X,r)$ is a \emph{cyclic permutation solution}. A permutation solution is indecomposable if and only if it is cyclic \cite[p.184]{etingof}.
\begin{prop}\cite{etingof}\label{prop-etingof-indecomposable}
	Let  $(X,r)$ be a non-degenerate involutive set-theoretic solution of the YBE, with structure group $G(X,r)$.   Then $(X,r)$  is indecomposable if and only if  $G(X,r)$ acts transitively on $X$.
\end{prop}
\begin{defn}
	The \emph{retract} relation $\sim$ on the set $X$ is defined by $x \sim y$ if $\sigma_x=\sigma_y$. There is a natural induced solution $Ret(X,r)=(X/\sim,r)$, called the \emph{the retraction of $(X,r)$}, defined by $r'([x],[y])=([\sigma_{x}(y)],[\gamma_y(x)])$.
	A non-degenerate involutive set-theoretic solution  $(X,r)$   is called \emph{a multipermutation solution of level $\ell$} if $\ell$ is the smallest natural number such that the solution $\mid \operatorname{Ret}^\ell(X,r)\mid=1$, where 
	$\operatorname{Ret}^k(X,r)=\operatorname{Ret}(\operatorname{Ret}^{k-1}(X,r))$, for $k>1$. If such an $\ell$ exists, $(X,r)$   is also called \emph{retractable}, otherwise it is called \emph{irretractable}.
\end{defn}
An important characterisation of  non-degenerate involutive set-theoretic solutions of the YBE is  presented in the following proposition.
\begin{thm}\cite[p.176-180]{etingof}\label{prop-etingof}
Let  $(X,r)$ be a non-degenerate involutive set-theoretic solution of the YBE, defined by  $r(x,y)=(\sigma_{x}(y),\gamma_{y}(x))$,   $x,y \in X$,  with structure group $G(X,r)$. Let $\mathbb{Z}^{X}$ denote the free abelian group with basis $\{t_x \mid x \in X\}$, and  $\operatorname{Sym}_X$ denote the symmetric group of $X$. Then 
\begin{enumerate}[(i)]
		\item The map $\varphi: G(X,r) \rightarrow \operatorname{Sym}_X$, defined by $x \mapsto \sigma_{x}$,  is a homomorphism of groups.
		\item  The group $\operatorname{Sym}_X$  acts on $\mathbb{Z}^{X}$.
		\item The group $G(X,r) $  acts on $\mathbb{Z}^{X}$:  if  $g \in G$, then $g \bullet t_x=t_{\alpha(x)}$, with $\alpha=\varphi(g)$. 
	\item The map 	$\pi:G(X,r)\rightarrow \mathbb{Z}^{X}$ is a bijective $1$-cocycle, where  $\pi(x)=t_x$,  for $x \in X$, and  $\pi(gh) =\pi(g)+g \bullet \pi(h)$,  for $g,h \in G(X,r)$.
	\item There  is a monomorphism of groups $\psi: G(X,r) \rightarrow \mathbb{Z}^{X} \rtimes \operatorname{Sym}_X$:  $\psi(x)=(t_x,\sigma_x) $, $\psi(g)=(\pi(g), \varphi(g))$.
 \end{enumerate}
\end{thm}
\begin{lem}\label{lem-braided--}
	Let $(X,r)$ be a non-degenerate involutive  set-theoretic solution. Then 
	\begin{equation}\label{eqn-braided-g--}
	\sigma_i\,\sigma_{\sigma_i^{-1}(j)}\,=\,\sigma_j\,\sigma_{\sigma_j^{-1}(i)} 
	\end{equation}
\end{lem}
\begin{proof}
	From the definition of $r$,  $r(x_i,x_{\sigma_i^{-1}(j)})=(x_{\sigma_i\sigma_i^{-1}(j)}, x_{\gamma_{\sigma_i^{-1}(j)}(i)})=(x_j, x_{\gamma_{\sigma_i^{-1}(j)}(i)})$. By replacing $x$ by  $i$ and $y$ by  $\sigma_i^{-1}(j)$  in Equation (\ref{eqn-inv1}), we have  $\gamma_{\sigma_i^{-1}(j)}(i)=\sigma_j^{-1}(i)$, that is $r(x_i,x_{\sigma_i^{-1}(j)})=(x_j, x_{\sigma_j^{-1}(i)})$.  So, the equality  $x_i\,x_{\sigma_i^{-1}(j)}=x_j\, x_{\sigma_j^{-1}(i)}$ holds in the structure group of $(X,r)$.  From Theorem  \ref{prop-etingof} $(i)$, this implies  $\sigma_i\,\sigma_{\sigma_i^{-1}(j)}\,=\,\sigma_j\,\sigma_{\sigma_j^{-1}(i)}$.
\end{proof}

\begin{lem}\cite{etingof}\label{lem_etingof_inverseT}
	Let $D: X \rightarrow X$ be the map defined by $D(x) =\sigma^{-1}_x(x)$. Then the  map $D$ is invertible and  $D^{-1}(y)= \gamma^{-1}_y(y)$,  with $x,y \in X$ such that   $r(x,y)=(x,y)$.
	Inductively, $D^{m}(x)=\sigma^{-1}_{D^{m-1}(x)}D^{m-1}(x)=\sigma^{-1}_{D^{m-1}(x)}\sigma^{-1}_{D^{m-2}(x)}..\sigma^{-1}_{D(x)}\sigma^{-1}_{x}(x)$.
\end{lem}
\begin{defn}\label{defn_class}
	$(i)$ We say that $(X,r)$  \emph{satisfies $(C)$}, if $\sigma_x\sigma_y=Id_X$, whenever $r(x,y)=(x,y)$;  $xy$ and $yx$ are called \emph{frozen elements of length $2$}  \cite{chou_godel2}.\\
	$(ii)$ We say that $(X,r)$ is of \emph{class $m$}, if $m$  is the  minimal  natural number such that \\
	$\sigma_x\,\sigma_{D(x)}\,\sigma_{D^{2}(x)}\,...\,\sigma_{D^{m-1}(x)}\,=\,Id_X$,  for every $x \in X$.
\end{defn}

\begin{rem}\label{remark_defn_class}
	\begin{enumerate}[(i)]
		\item 	In  Definition \ref{defn_class}$(ii)$, we use the terminology from   \cite{deh_coxeterlike}, with a different (but equivalent) formulation. Being of class $2$ is equivalent to satisfying  $(C)$, with $(C)$ rewritten as  $\sigma_x\sigma_{\sigma_x^{-1}(x)}=Id_X$,  $\forall x \in X$.

\item For a non-degenerate and involutive  set-theoretic solution $(X,r)$ of class $m$,  with $\mid X\mid=n$,  there exists  a subgroup $N$ of $G(X,r)$, generated by the $n$ frozen elements of length $m$,   which  is  normal, free abelian  of rank $n$ and the group $W$ defined by $G(X,r)/N$ is finite of order $m^n$ (and is a  Coxeter-like  group)  \cite{deh_coxeterlike} (and \cite{chou_godel2} for the case $m=2$).  It is also proved there  that $N \subset \operatorname{Ker}(\varphi)$, where $\varphi: G(X,r) \rightarrow \operatorname{Sym}_X$ is the homomorphism defined by $\varphi(x)=\sigma_{x}$  in  Theorem \ref{prop-etingof}$(i)$. So, 
the IYB group of $(X,r)$, $\mathcal{G}(X,r)$,  is a subgroup of  the group $W$,  and  its order divides $m^n$ from Lagrange's Theorem.
	\end{enumerate}
\end{rem}
\begin{ex} \label{exemple-irretractable+indecomposable4} Let $X = \{x_1,x_2,x_3,x_4\}$, and ~$r: X\times X\to X\times X$ be defined by $r(x_i,x_j)=(x_{\sigma_{i}(j)},x_{\gamma_{j}(i)})$,  where~$\sigma_i$ and $\gamma_j$ are permutations on~$\{1,2,3,4\}$ as follows: $\sigma_{1}=(1,2,3,4)$, $\sigma_2=(2,1,4,3)$, $\sigma_{3}=(1,3)(2)(4)$, $\sigma_{4}=(1)(3)(2,4)$; $\gamma_{1}=(1,2,4,3)$, $\gamma_2=(2,1,3,4)$, $\gamma_{3}=(1)(2,3)(4)$, $\gamma_{4}=(1,4)(2)(3)$.
	Then  $(X,r)$ is  an indecomposable, and irretractable  solution, with   structure group $G=\operatorname{Gp} \langle X\mid  
	x_{1}^2=x_{2}^2;
x_{1}x_{2}=	x^{2}_{3};
 x_{1}x_{3}=x_{4}x_{1};
	x_{2}x_{1}=x^{2}_{4}; 
	x_{2}x_{4}=x_{3}x_{2};
	x_{3}x_{4}=x_{4}x_{3} \rangle$.  The four frozen elements of length $2$ are  $x_{1}x_{4}$, $x_{4}x_{2}$, $ x_{2}x_{3}$ and $x_{3}x_{1}$. The class of the solution is $4$ since $\sigma_{1}\sigma_{4}\sigma_{2}\sigma_{3}=Id_X$, and  the  four frozen elements of length $4$ are  $x_{1}x_{4}x_{2}x_{3}$,  $x_{4}x_{2}x_{3}x_{1}$,   $x_{2}x_{3}x_{1}x_{4}$, and   $x_{3}x_{1}x_{4}x_{2}$. The quotient group $W$ has order $4^4$ and the IYB group $\mathcal{G}(X,r)$  has order $8$.
\end{ex}    
\subsection{The Yang-Baxter equation and the Faddeev-Reshetikhin-Takhtadjian construction}

In \cite{drinf-quasi1,drinf-quasi2}, Drinfeld introduced the concept of  braided (or quasi-triangular) bialgebras. These are bialgebras with a  universal $R$-matrix inducing a solution of  the  Yang-Baxter equation on any of their modules. There is a dual notion of cobraided bialgebras, these are bialgebras inducing 
$R$-matrices on their comodules.   Both constructions provide a systematic method to produce solutions of   the    Yang-Baxter equation. Conversely, given a  solution of  the Yang-Baxter equation, 
 there is a method, due to Faddeev, Reshetikhin and  Takhtadjian (FRT construction for short),  to construct a cobraided bialgebra  \cite{frt}, \cite{kassel}. Faddeev, Reshetikhin and Takhtadjian show that if  $c \in \operatorname{End}(V \otimes V)$, where  $V$  is a finite-dimensional vector space, is a solution of the  Yang-Baxter equation, then there exists  a cobraided bialgebra $A(c)$ coacting on $V$ \cite{frt}.  We describe some of the main steps of the FRT construction  and  refer the reader to \cite{kassel} [VIII] for a complete exposition of the FRT construction and for more details in the topic.\\

Let $\{v_i \mid 1 \leq i \leq n\}$ be a basis of $V$ and let   the coefficients $ c_{ij}^{\ell m} $ be defined by
\begin{equation}\label{eqn-coeff}
c(v_i \otimes v_j)= \sum\limits_{1 \leq \ell,m \leq n} c_{ij}^{\ell m} v_{\ell} \otimes v_m
\end{equation}

\begin{defn}\cite{frt}, \cite{kassel}[VIII.6]
The algebra  $A(c)$ is the quotient of the free algebra $F$ generated by the set  $\{T_i^k \mid1 \leq i,k \leq n \} $ by the two-sided  ideal $I(c)$ generated by all the elements $C_{ij}^{kl}$, where 
\begin{equation}\label{eqn-gen-ideal}
C_{ij}^{kl}= \sum\limits_{1 \leq \ell,m \leq n}     c_{ij}^{\ell m}T_\ell^k T_m^l\;
-\; 
\sum\limits_{1 \leq \ell,m \leq n}    T_i^\ell T_j^m   c_{\ell m}^{kl}
\end{equation}
	and $i,j,k,l$ run over the indexing set.
\end{defn}
\begin{thm}\cite{frt},\cite{kassel}[VIII.6]
Let $V$ be a vector space of dimension $n$. Let $c \in \operatorname{End}(V \otimes V)$ a solution of the Yang-Baxter equation.  Then   
\begin{enumerate}[(i)]
	\item there exists a  unique bialgebra structure on $A(c)$, such that\\ $\Delta(T_i^k)=\sum\limits_{1 \leq \ell \leq n}T_i^\ell \otimes T_\ell^k$ and $\epsilon(T_i^k)=\delta_{ik}$.
	\item  there exists  a linear map $\Delta_V:V \rightarrow A(c)\otimes V$ such that   $\Delta_V$ equips $V$ with the structure of a comodule over $A(c)$ and the map $c$ becomes a comodule map with respect to this structure.
	\item there exists a unique linear form $r$ on  $A(c)\otimes A(c)$, such that $A(c)$ is a cobraided bialgebra and $r(T_i^k\otimes T_j^l)=c_{ji}^{kl}$, for every $ 1 \leq i,j,k,l \leq n$.
\end{enumerate}
\end{thm}
\section{Some preliminaries on Cryptography}
We use the terminology from \cite{diffie} and \cite{rsa}, and   also from the books \cite{book-crypto}, \cite{book2crypto}.  The literature in this topic is very vast, and the interested reader can find many expositions of the topic. In \cite{diffie},  Diffie and Hellman  introduced  the concept of  a “public-key cryptosystem”, yet without any practical implementation of such a system.  In \cite{rsa}, Rivest, Shamir and Adelman presented  a public-key cryptosystem, that is they constructed  an encryption method  with the property that publicly revealing an encryption key does not reveal the corresponding decryption key.  From that time, this field of research has very much developed and there are several crytographical problems addressed. In the context of  Public key Encryption, the following definition is often used.
\begin{defn}\cite{diffie},\cite{rsa}
	Let $E$   and $D$ denote the encryption  and decryption procedure respectively.
	The function $E$ is called a \emph{trap-door one-way function} if $E$ satisfies the following properties:\begin{enumerate}[(i)]
		\item $D(E(M) =M$, i.e. deciphering the enciphered form of a message $M$ yields $M$. 
		\item Both $E$ and $D$ are easy to compute. 
		\item  By publicly revealing $E$,  the user does not reveal an easy way to compute $D$. This means that in practice only he can decrypt messages encrypted with $E$, or compute $D$ efficiently.
	\end{enumerate}
	A trap-door one-way function is called  a \emph{trap-door one-way permutation} if 
	it satisfies $E(D(M)) = M$. That is, if a message $M$ is first deciphered and then enciphered, $M$ is the result. 
\end{defn}
The idea is, that  within the mathematics, it should be very difficult to determine the private key given the public key, such as factorising a very large number into its prime factors. Three main  methods used for this include integer factorisation (such as RSA \cite{rsa}), discrete logarithms (as ElGamal \cite{elgamal}), and elliptic curve relationships (as the elliptic curve \cite{elliptic1,elliptic2}).  The four main problems in cryptography are: encryption/decryption, key exchange, authentication and signature. The following idea is common to all the problems: there are two entities, traditionally called A(lice) and B(ob),  who want to communicate  in such a way that an intruder observing the communication could not understand.  We describe each of them  briefly.\\

\textbf{ Public key Encryption:}
It is an asymmetric key method, as it uses a public key (which can be distributed) and a private key (which should be kept secret).  Bob wishes to send Alice a message $M$, and he can use  the public key to encrypt  his message. Alice must be able to retrieve Bob’s original message using her private key, but an intruder watching the communication should not.\\

\textbf{Key Exchange:}
Alice and Bob  wish to agree on a common secret, in such a way that an intruder observing the communication cannot deduce any useful information about the common secret.\\
\textbf{Authentication:}
Alice (the prover) wishes to prove her identity to Bob (the verifier), i.e., she wishes to prove that she knows some private (secret) key without enabling an intruder watching the communication to deduce anything about her private key.\\
\textbf{Signature:}
Alice  wishes to  send  Bob a (clear or ciphered) message together with a signature proving the origin of the message.  To implement signatures the public-key cryptosystem must be implemented with trap-door one-way permutations. How can user Bob send Alice a “signed” message M in a public-key cryptosystem?  The following procedure is described in \cite{rsa} in the following way.
Let  $D_A$, $E_A$, $D_B$ and  $E_B$,  denote the decryption and encryption procedures of Alice  and Bob respectively.  Bob first computes his “signature” $S$ for the message $M$ using $D_B$.  He computes: $S= D_B(M)$.
He then encrypts $S$ using $E_A$, and sends the result $E_A(S)$ to Alice. He need not send $M$, as  it can be computed from $S$.
Alice first decrypts the ciphertext with $D_A$ to obtain $S$. She knows who is the presumed sender of the signature (in this case, Bob); this can be given if necessary in plain text attached to $S$. She then extracts the message with the encryption procedure of the sender, in this case $E_B$, and computes $M = E_B(S)$.

Therefore Alice has received a message “signed” by Bob, which she can “prove” that he sent, but which she cannot modify.

\section{Iterative construction of set-theoretic solutions and proof of Theorem \ref{thmA-construction} }

\subsection{Application of the FRT construction on a set-theoretic solution}
\label{subsec-frt-set-theoretic}
Let $(X,r)$ be a non-degenerate and involutive  solution of the  Yang-Baxter equation,  where $\mid X\mid =n$ and $r(x,y)=(\sigma_{x}(y),\gamma_{y}(x))$, $x,y \in X$.  Let $c \in \operatorname{End}(V \otimes V)$ be  induced by  $(X,r)$,   where  $V$ is  the vector space spanned by $X$.  It is interesting to understand how the FRT construction applies in this case. This question is addressed  in \cite{etingof}, where the authors compute the Hilbert series of the algebra $A(c)$. In \cite{doikou}, the authors study the  quantum groups associated with Baxterized solutions of the Yang–Baxter equation coming from braces,  via the FRT construction. Here, we are interested in a more combinatorial approach  and in particular we are interested  to find  how the generators of the two-sided ideal $I(c)$ look like.
\begin{lem}\label{lem-frt-set-theo}
	The algebra  $A(c)$ is the quotient of the free algebra $F$ generated by the set  $\{T_i^k \mid1 \leq i,k \leq n \} $ by the two-sided  ideal $I(c)$
	generated by all the elements $C_{ij}^{kl}$, where 
	\begin{equation}\label{eqn-ideal-set-theoretic}
	C_{ij}^{kl}\,=\, T_{\sigma_i(j)}^{k}T_{\gamma_j(i)}^{l}\,-\,T_{i}^{\sigma_k(l)}T_{j}^{\gamma_l(k)}
	\end{equation}
	\begin{equation}\label{eqn-ideal-set-theoretic2}
	C_{\sigma_i(j)\gamma_j(i)}^{kl}\,=\, T_{i}^{k}T_{j}^{l}\,-\,T_{\sigma_i(j)}^{\sigma_k(l)}T_{\gamma_j(i)}^{\gamma_l(k)}
	\end{equation}
	Furthermore, the following hold:
	\begin{enumerate}[(i)]
		\item  $  C_{\sigma_i(j)\gamma_j(i)}^{\sigma_k(l)\gamma_l(k)}\,\,\,=\,-\, C_{ij}^{kl} $.
		
		\item $  C_{ij}^{\sigma_k(l)\gamma_l(k)}\,\,\,=\,-\, C_{\sigma_i(j)\gamma_j(i)}^{kl}$.
		
		\item if $r(i,j)=(i,j)$, then 
		$  C_{ij}^{\sigma_k(l)\gamma_l(k)}\,\,\,=\,\pm\, C_{ij}^{kl} $.
		\item if $r(k,l)=(k,l)$, then $  C_{\sigma_i(j)\gamma_j(i)}^{kl}\,\,\,=\,\pm\, C_{ij}^{kl} $.
		\item if $r(i,j)=(i,j)$ and $r(k,l)=(k,l)$, then $C_{ij}^{kl}=0 $.
		\item  there are $n^2\choose 2$ elements in the set $\{C_{ij}^{kl}\}/\pm$.
	\end{enumerate}
\end{lem}
\begin{proof}
	In order to  compute $C_{ij}^{kl}$, we need to compute the coefficients $c_{ij}^{\ell m}$ as defined in Equation (\ref{eqn-coeff}). For a set-theoretic solution,  we have $r(x_i,x_j)=(x_{\sigma_{i}(j)},x_{\gamma_{j}(i)})$,  so $c(v_i \otimes v_j)=v_{\sigma_{i}(j)}\,\otimes\,v_{\gamma_{j}(i)}$, that is $c_{ij}^{\sigma_{i}(j)\gamma_{j}(i)}=1$ and $c_{ij}^{\ell m}=0$,  if $\ell \neq \sigma_{i}(j)$ or $m \neq \gamma_{j}(i)$. From the involutivity of the solution, we have $c_{\sigma_{i}(j)\gamma_{j}(i)}^{ij}=1$ and $c_{\sigma_{i}(j)\gamma_{j}(i)}^{\ell m}=0$,  if $\ell \neq i$ or $m \neq j$.  So, using that in  Equation (\ref{eqn-gen-ideal}), we have Equations (\ref{eqn-ideal-set-theoretic}) and (\ref{eqn-ideal-set-theoretic2}).\\
	$(i)$ and  $(ii)$  hold  from the involutivity of the solution. Indeed, 
	$(i)$ results from substituting $\sigma_k(l)\gamma_l(k)$ instead of $kl$ in Equation (\ref{eqn-ideal-set-theoretic2}) and using Equations (\ref{eqn-inv1})-(\ref{eqn-inv2}) and	$(ii)$  results from substituting $\sigma_k(l)\gamma_l(k)$ instead of $kl$ in Equation (\ref{eqn-ideal-set-theoretic}) and using  (\ref{eqn-inv1})-(\ref{eqn-inv2}). \\
	$(iii)$, $(iv)$, $(v)$  result from  $(i)$-$(ii)$  and the fact that if   $r(i,j)=(i,j)$, then $\sigma_i(j)\gamma_j(i)=ij$.\\
	$(vi)$ From  Equation (\ref{eqn-ideal-set-theoretic2}), whenever  the pair $(T_{i}^{k},T_{j}^{l})$ is chosen,  the pair $(T_{\sigma_i(j)}^{\sigma_k(l)},T_{\gamma_j(i)}^{\gamma_l(k)})$ is uniquely  determined. Let renumber the set $\{T_i^k \mid1 \leq i,k \leq n \} $ in the following order: $T_1^1,T_1^2,...,T_1^n,T_2^1,....,T_n^1,...,T_n^n$, that is $T_1^1$ is the first element,  $T_2^1$ is the $n+1$-th element, $T_n^n$ is the $n^2$-th element and so on. So, the number of elements in the set $\{C_{ij}^{kl}\}/\pm$ is equal to the number of possibilities to choose a  pair $(T_{i}^{k},T_{j}^{l})$  from the set $T_1^1,T_1^2,...,T_1^n,T_2^1,....,T_n^1,...,T_n^n$, that is  $n^2\choose 2$.
\end{proof}
\subsection{Construction of non-degenerate and involutive  set-theoretic solutions}\label{subsec-iterative}
  In Section \ref{subsec-frt-set-theoretic}, we computed the generating elements, $	C_{ij}^{kl}$, of the two-sided ideal $I(c)$ in the cobraided bialgebra $A(c)$, whenever $c$ is a solution of the YBE induced by  a non-degenerate and involutive solution $(X,r)$, with $\mid X\mid=n$. We now forget about the  algebraic structures  $A(c)$ and  $I(c)$, and consider only the equations derived from them. In particular, we consider Equation (\ref{eqn-ideal-set-theoretic2}), and define from it a new  set-theoretic solution of size $n^2$.
 \begin{defn} \label{defn-new-solution}
 	Let $(X,r)$ be a set-theoretic solution, with $\mid X\mid=n$. Let 
 	$\mathcal{X}^2$ denote the set of elements  $\{T_{i}^{k}\mid 1 \leq i,k \leq n\}$ 
 	in bijection with $X \times X$.  We define the following map:
 	\begin{gather}
 \tilde{r}: \mathcal{X}^2  \times \mathcal{X}^2  \rightarrow\mathcal{X}^2 \times \mathcal{X}^2 \nonumber\\
  \tilde{r}(T_{i}^{k},T_{j}^{l})\,=\,(T_{\sigma_i(j)}^{\sigma_k(l)}\,,\,T_{\gamma_j(i)}^{\gamma_l(k)}) \label{eqn-defn-new-r}
 	\end{gather}
 	 We write $ \tilde{r}(T_{i}^{k},T_{j}^{l})\,=(g_i^k(T_{j}^{l})\,,\,f_j^l(T_{i}^{k}))$,  where   $g_i^k$ and $f_j^l$, $1 \leq i,j,k,l \leq n$ are defined by:
 	 \begin{gather}
 	g_i^k\,,\, f_j^l\,: \mathcal{X}^2  \rightarrow\mathcal{X}^2  \nonumber\\
 	 g_i^k(T_{j}^{l})\,=\, T_{\sigma_i(j)}^{\sigma_k(l)} \label{eqn-defn-new-g}\\
 	 f_j^l(T_{i}^{k})\,=\,T_{\gamma_j(i)}^{\gamma_l(k)} \label{eqn-defn-new-f}
 	  	 \end{gather}
 	 \end{defn}
Note that if $(X,r)$ is the trivial solution of size $n$, $(\mathcal{X}^2,  \tilde{r})$ 
is the trivial solution of size $n^2$.
 \begin{ex}\label{ex-permut-sol-2}
 	Let $(X,r)$ be  an indecomposable permutation solution, with $X=\{x_1,x_2\}$, and  $\sigma_1=\sigma_2=\gamma_1=\gamma_2=(1,2)$. Its structure group is $G(X,r)=\langle x_1,x_2 \mid x_1^2=x_2^2\rangle $.  Then,  $\mathcal{X}^2=\{T_1^1,T_1^2,T_2^1,T_2^2\}$  and $g_i^k=f_j^l=\,(T_1^1,T_2^2)(T_1^2,T_2^1)$, for  $1 \leq i,j,k,l \leq 2$.
 \end{ex}
 In the following, we show that the pair $(\mathcal{X}^2,  \tilde{r})$ satisfies the same properties as $(X,r)$, that is  if  $(X,r)$ is  non-degenerate and involutive, with $\mid X\mid=n$, then  $(\mathcal{X}^2,  \tilde{r})$  is also  non-degenerate and involutive,  with $\mid \mathcal{X}^2\mid=n^2$.
 \begin{lem}\label{lem-new--invol+braided}
 	Let $(X,r)$ be a non-degenerate and involutive set-theoretic solution, with $\mid X\mid=n$. Let  the pair $(\mathcal{X}^2,  \tilde{r})$,  and the functions 
 $	g_i^k\,,\, f_j^l\,: \mathcal{X}^2  \rightarrow\mathcal{X}^2$ be defined as in Definition \ref{defn-new-solution}. Then 
 	\begin{enumerate}[(i)]
 		\item  $ \tilde{r}$  is bijective.
 	\item  $(\mathcal{X}^2,  \tilde{r})$ is non-degenerate, that is $g_i^k\,,\, f_j^l$, $1 \leq i,j,k,l \leq n$, are bijective.
 			\item   $(\mathcal{X}^2,  \tilde{r})$ is  involutive, that is  $\tilde{r}^2=\,Id_{\mathcal{X}^2\times \mathcal{X}^2}$. 			
 	\item    $(\mathcal{X}^2,  \tilde{r})$ is  braided, that is,  $\tilde{r}^{12}\tilde{r}^{23}\tilde{r}^{12}\,=\,\tilde{r}^{23}\tilde{r}^{12}\tilde{r}^{23}$. 	Additionally, for every  $1 \leq i,j,k,l,s,m \leq n$,   the following  equations hold:\\
 	 \begin{gather}
 	g_i^kg_j^l\,=\, g_{\sigma_i(j)}^{\sigma_k(l)}\,g_{\gamma_j(i)}^{\gamma_l(k)} \label{eqn-braided-g}\\
 		f_j^lf_i^k\,=\, f_{\gamma_j(i)}^{\gamma_l(k)} \,f_{\sigma_i(j)}^{\sigma_k(l)} \label{eqn-braided-f}\\
 			f_{\sigma_{\gamma_j(i)}(s)}^{\sigma_{\gamma_l(k)}(m)}\,g_i^k(T_j^l)\;=\;
 		g_{\gamma_{\sigma_j(s)}(i)}^{\gamma_{\sigma_l(m)}(k)}\,\,f_s^m(T_j^l)
 		\label{eqn-braided-new}
 	\end{gather} 
 	\end{enumerate}
 \end{lem}
\begin{proof}
	$(i)$, $(ii)$   From the definition of  $g_i^k$, $g_i^k(T_{j}^{l})\,=\, T_{\sigma_i(j)}^{\sigma_k(l)}$ and $g_i^k(T_{s}^{m})\,=\, T_{\sigma_i(s)}^{\sigma_k(m)}$. As $(X,r)$ is non-degenerate,  $\sigma_i$, $\sigma_k$ are bijective, and 
	 $ T_{\sigma_i(j)}^{\sigma_k(l)}\neq  T_{\sigma_i(s)}^{\sigma_k(m)}$, if $j \neq s$ or $l\neq m$.  The bijectivity of   the functions $f_j^l$ relies on the bijectivity of $\gamma_j$ and $\gamma_l$. So, $(\mathcal{X}^2,  \tilde{r})$ is non-degenerate 
	and  $ \tilde{r}$  is bijective.\\
	$(iii)$ $\tilde{r}^2(T_{i}^{k},T_{j}^{l})\,= \tilde{r}\,(T_{\sigma_i(j)}^{\sigma_k(l)}\,,\,T_{\gamma_j(i)}^{\gamma_l(k)}) \,=\, (T_{\sigma_{\sigma_i(j)}\gamma_j(i)}^{\sigma_{\sigma_k(l)}\gamma_l(k)}\,,\, 
	T_{\gamma_{\gamma_j(i)}\sigma_i(j)}^{\gamma_{\gamma_l(k)}\sigma_k(l)})$. As $(X,r)$ is involutive,  this is equal to $(T_{i}^{k},T_{j}^{l})$,  from Equations (\ref{eqn-inv1}), (\ref{eqn-inv2}). So,     $\tilde{r}^2=\,Id_{\mathcal{X}^2\times \mathcal{X}^2}$. 			\\
	$(iv)$  From the definition of $\tilde{r}$, 	$\tilde{r}^{12}\tilde{r}^{23}\tilde{r}^{12}\,=\,\tilde{r}^{23}\tilde{r}^{12}\tilde{r}^{23}$ if and only if Equations (\ref{eqn-braided-g})-(\ref{eqn-braided-new}) hold. We prove (\ref{eqn-braided-g}). From  Equation (\ref{eqn-defn-new-g}), we have:
		\begin{align*}
	g_i^kg_j^l(T_s^m)\,=\,g_i^k(T_{\sigma_j(s)}^{\sigma_l(m)})=\,
	T_{\sigma_i\sigma_j(s)}^{\sigma_k\sigma_l(m)}\\
	g_{\sigma_i(j)}^{\sigma_k(l)}\,g_{\gamma_j(i)}^{\gamma_l(k)}(T_s^m)=
	g_{\sigma_i(j)}^{\sigma_k(l)}(T_{\sigma_{\gamma_j(i)}(s)}^{\sigma_{\gamma_l(k)}(m)})\,=\,T_{\sigma_{\sigma_i(j)}\sigma_{\gamma_j(i)}(s)}^{\sigma_{\sigma_k(l)}\sigma_{\gamma_l(k)}(m)}
		\end{align*}
	From Equation (\ref{eqn-braided-sigma}), 
	$T_{\sigma_i\sigma_j(s)}^{\sigma_k\sigma_l(m)}\,=\,T_{\sigma_{\sigma_i(j)}\sigma_{\gamma_j(i)}(s)}^{\sigma_{\sigma_k(l)}\sigma_{\gamma_l(k)}(m)}$,  for every $1 \leq s,m \leq n$, so (\ref{eqn-braided-g}) holds. In the same way, we show (\ref{eqn-braided-f}) holds, using Equations   (\ref{eqn-defn-new-f}) and (\ref{eqn-braided-gamma}).  We prove (\ref{eqn-braided-new}):		
	\begin{align*}
f_{\sigma_{\gamma_j(i)}(s)}^{\sigma_{\gamma_l(k)}(m)}\,g_i^k(T_j^l)\;=
f_{\sigma_{\gamma_j(i)}(s)}^{\sigma_{\gamma_l(k)}(m)}\,(T_{\sigma_i(j)}^{\sigma_k(l)})\;=T_{\gamma_{\sigma_{\gamma_j(i)}(s)}\sigma_i(j)}^{\gamma_{\sigma_{\gamma_l(k)}(m)}\sigma_k(l)}\\
g_{\gamma_{\sigma_j(s)}(i)}^{\gamma_{\sigma_l(m)}(k)}\,\,f_s^m(T_j^l)\,	=\,g_{\gamma_{\sigma_j(s)}(i)}^{\gamma_{\sigma_l(m)}(k)}\,\,(T_{\gamma_s(j)}^{\gamma_m(l)})\,=\,	T_{\sigma_{\gamma_{\sigma_j(s)}(i)} \gamma_s(j)}^{\sigma_{\gamma_{\sigma_l(m)}(k)} \gamma_m(l)}
	\end{align*}
From Equation (\ref{eqn-braided-old}), these are equal, that is (\ref{eqn-braided-new})
 holds.				
\end{proof}

\begin{lem}\label{lem-new-class}
	Let $(X,r)$ be a non-degenerate and involutive set-theoretic solution, with $\mid X\mid=n$. Let  the pair $(\mathcal{X}^2,  \tilde{r})$,  and the functions 
$	g_i^k\,,\, f_j^l\,: \mathcal{X}^2  \rightarrow\mathcal{X}^2$ be defined as in Definition \ref{defn-new-solution}. If $(X,r)$  is of class $m$, $m>1$, then $(\mathcal{X}^2,  \tilde{r})$ is also of class $m$.
\end{lem}
\begin{proof}
We  recall that $(X,r)$ is of  class $m$, if $m$  is the  minimal  natural number such that 
$\sigma_x\,\sigma_{D(x)}\,\sigma_{D^{2}(x)}\,...\,\sigma_{D^{m-1}(x)}\,=\,Id_X$,  for every $x \in X$, with  $D(x) =\sigma^{-1}_x(x)$.  
From the definition of $\tilde{r}$, $\tilde{r}(T_x^y,\,T_{\sigma^{-1}_x(x)}^{\sigma^{-1}_y(y)})=(T_x^y,\,T_{\sigma^{-1}_x(x)}^{\sigma^{-1}_y(y)})$, that is   $T_x^yT_{\sigma^{-1}_x(x)}^{\sigma^{-1}_y(y)}\,=\,T_x^yT_{D(x)}^{D(y)}$ is a frozen element of length $2$.  In the same way,  $\tilde{r}(T_{D(x)}^{D(y)},\, T_{D^2(x)}^{D^2(y)})\,=\,(T_{D(x)}^{D(y)},\, T_{D^2(x)}^{D^2(y)})$,  and so on  $\tilde{r}(T_{D^{k}(x)}^{D^{k}(y)},\, T_{D^{k+1}(x)}^{D^{k+1}(y)})\,=\,(T_{D^{k}(x)}^{D^{k}(y)},\, T_{D^{k+1}(x)}^{D^{k+1}(y)})$.  So, we show that, for every $x,y\in X$ ($x,y$ not necessarily  distinct),  $m$  is the  minimal  natural number such that 
$g_x^y\,g_{D(x)}^{D(y)}\,g_{D^{2}(x)}^{D^{2}(y)}\,...\,g_{D^{m-1}(x)}^{D^{m-1}(y)}\,=\,Id_{\mathcal{X}^2}$.  Let $T_i^k \in \mathcal{X}^2$. Then,  from   Definition \ref{defn-new-solution}:\\
$g_x^y\,g_{D(x)}^{D(y)}\,g_{D^{2}(x)}^{D^{2}(y)}\,...\,g_{D^{m-1}(x)}^{D^{m-1}(y)}(T_i^k)\;=\;T_{\sigma_x\,\sigma_{D(x)}\,\sigma_{D^{2}(x)}\,...\,\sigma_{D^{m-1}(x)}(i)}^{\sigma_y\,\sigma_{D(y)}\,\sigma_{D^{2}(y)}\,...\,\sigma_{D^{m-1}(y)}(k)}\;=\; T_i^k $, since  $(X,r)$ is of  class $m$, and $m$ is the minimal such number. That is, $(\mathcal{X}^2,  \tilde{r})$ is also of class $m$.
	\end{proof}
From the proof of Lemma \ref{lem-new-class}, for every $x,y\in X$ ($x,y$ not necessarily  distinct), the frozen elements of length $m$ have the form $T_x^y\,T_{D(x)}^{D(y)}\,T_{D^{2}(x)}^{D^{2}(y)}\,...\,T_{D^{m-1}(x)}^{D^{m-1}(y)}$.
 \begin{ex}\label{ex-continued-n=2}
The pair $(\mathcal{X}^2,  \tilde{r})$  from Example \ref{ex-permut-sol-2} is a non-degenerate and involutive solution of size 4. It is decomposable with 
$\mathcal{X}^2=\{T_1^1,T_2^2\}\cup \{T_1^2,T_2^1\}$.  The class of the solution is $2$ and the four frozen elements of length $2$ are  $T_1^1T_2^2$, $T_2^2T_1^1$, $T_2^1T_1^2$ and $T_1^2T_2^1$. The  structure group 
 $G(\mathcal{X}^2,  \tilde{r})= \langle T_1^1,T_1^2,T_2^1,T_2^2\mid  (T_1^1)^2=(T_2^2)^2, 
 T_1^1T_1^2=T_2^1T_2^2,
T_1^2T_1^1=T_2^2T_2^1, 
T_1^1T_2^1=T_1^2T_2^2,
(T_1^2)^2=(T_2^1)^2,
T_2^1T_1^1=T_2^2T_1^2\rangle$.  
 \end{ex}
\begin{lem}\label{lem-equiv}
Let $(X,r)$ and $(X',r')$ be  non-degenerate and involutive set-theoretic solutions of the   Yang-Baxter equation, with $\mid X\mid=\mid X'\mid=n$.
\begin{enumerate}[(i)]
	\item If $(X,r)$ and $(X',r')$  are different, then $(\mathcal{X}^2,  \tilde{r})$ and $(\mathcal{X'}^2,  \tilde{r}')$ are different.
	\item If $(X,r)$ and $(X',r')$  are isomorphic  then $(\mathcal{X}^2,  \tilde{r})$ and $(\mathcal{X'}^{2},  \tilde{r}')$ are isomorphic.
\end{enumerate} 
\end{lem}
\begin{proof}
$(i)$ results directly from Equation (\ref{eqn-defn-new-r}). Indeed, if $(X,r)$ and $(X',r')$  are different, then there exists $1 \leq i \leq n$ such that $\sigma_i$ and $\sigma'_i$ are different and so for every  $1 \leq k \leq n$, $g_{i}^{k}$ and $g_{i}^{'k}$ are different.\\
$(ii)$ Since  $(X,r)$ and $(X',r')$ are isomorphic, there is a bijection $\mu:X \rightarrow X'$ such that $(\mu \times \mu) \circ r=r'\circ (\mu\times \mu)$, that is for $x,y \in X$, if $\mu(x)=t$ and $\mu(y)=u$, then $\mu(\sigma_x(y))=\alpha_t(u) =\alpha_t(\mu(y))$  and $\mu(\gamma_y(x))=\beta_u(t) =\beta_u(\mu(x))$, where $\alpha_i,\beta_i$ are  the permutations defining $(X',r')$. We define $\tilde{\mu}:\mathcal{X}^2 \rightarrow\mathcal{X'}^2$ by  $\tilde{\mu}(T_{x_1}^{x_2})=S_{\mu(x_1)}^{\mu(x_2)}$, then a little computation shows that $\tilde{r}'\circ (\tilde{\mu}\times \tilde{\mu})\,=\,(\tilde{\mu}\times \tilde{\mu}) \circ\tilde{r}$, i.e.  $(\mathcal{X}^2,  \tilde{r})$ and $(\mathcal{X'}^2,  \tilde{r}')$ are isomorphic.
\end{proof}

 We give an easy method to compute the cycle decomposition of the bijection $g_i^k$ from the the cycle decompositions of  $\sigma_i$ and $\sigma_k$.  By  definition, $g_i^k(T_{j}^{l})\,=\, T_{\sigma_i(j)}^{\sigma_k(l)}$, so the length of the cycle with $T_{s_1}^{m_1}$ in  $g_i^k$ is the lcm of  the length of the cycle with $s_1$   in $\sigma_i$ and the length of the cycle  with $m_1$  in $\sigma_k$.  Assume that the  cycle in $\sigma_i$ with $s_1$ is $(s_1,s_2,..,s_q)$ and the cycle in $\sigma_k$ with $m_1$ is $(m_1,m_2,..,m_p)$.  Then, the cycle with $T_{s_1}^{m_1}$ in  $g_i^k$ has  the  form:\\
 
  \begin{align}
 (\;T_{s_1}^{m_1},T_{\sigma_i(s_1)}^{m_2},...,T_{\sigma_i^p(s_1)}^{m_p},T_{\sigma_i^{p+1}(s_1)}^{m_1},...,T_{\sigma_i^{q-1}(s_1)}^{m_p},...,T_{s_1}^{m_1},...,T_{\sigma_i^{q-1}(s_1)}^{m_p}\;) &&  \mathrm{if }\;\; q>p\\
 (\;T_{s_1}^{m_1},T_{s_2}^{\sigma_k(m_1)},...,T_{s_q}^{\sigma_k^{q}(m_1)},T_{s_1}^{\sigma_k^{q+1}(m_1)},...T_{s_q}^{\sigma_k^{p-1}(m_1)},...,T_{s_1}^{m_1},..., T_{s_q}^{\sigma_k^{p-1}(m_1)}\;)&& \mathrm{if }\;\; q< p
 \end{align}\\
If $q=p$, then  it has the form $ (\;T_{s_1}^{m_1},T_{s_2}^{m_2},.., T_{s_p}^{m_p}\;)$.  
\begin{rem}\label{rem-renumber}
The elements $T_i^k$  can be renumbered  using the following conversion rules:
\begin{gather*}
\textrm{Given }\, T_i^k ,\,\textrm{the corresponding  number is  }\,m=\,n(i-1)\,+\, k.\\
\textrm{Given }\, 1 \leq m \leq n^2, \,T_i^k  \, \textrm{is obtained in the following way:}\, i=\lceil\frac{m}{n} \rceil\; \textrm{and} \, k\equiv m\,(mod n).
\end{gather*}
 The permutation $g_i^k $ can be renumbered accordingly and it  can also be  rewritten  as a permutation in $S_{n^2}$.
\end{rem}
In the following example, we  illustrate the method of computation of  the functions $g_i^k$,  according to the 3 cases presented above. 
\begin{ex}\label{ex-pump4}
	For the solution in Example \ref{exemple-irretractable+indecomposable4}, we compute some of the functions $g_i^k$: 
	$g_1^3=(T_1^1,T_2^3,T_3^1,T_4^3)\;(T_1^3,T_2^1,T_3^3,T_4^1)\;(T_1^2,T_2^2,T_3^2,T_4^2)\;(T_1^4,T_2^4,T_3^4,T_4^4)$

$g_4^2=(T_1^1,T_1^4,T_1^3,T_1^2)\;(T_2^1,T_4^4,T_2^3,T_4^2)\;(T_3^1,T_3^4,T_3^3,T_3^2) \;(T_4^1,T_2^4,T_4^3,T_2^2)$

	$g_1^2=(T_1^1,T_2^4,T_3^3,T_4^2)\;(T_1^2,T_2^1,T_3^4,T_4^3)\;(T_1^3,T_2^2,T_3^1,T_4^4)\;(T_1^4,T_2^3,T_3^2,T_4^1)$
	
		$g_3^4=(T_1^1,T_3^1)\;(T_1^2,T_3^4)\;(T_1^3,T_3^3)\;(T_1^4,T_3^2)
		\;(T_2^1,T_2^1)	\;(T_2^2,T_2^4)	\;(T_2^3,T_2^3)\;(T_4^1,T_4^1)	\;(T_4^2,T_4^4)		\;(T_4^3,T_4^3)$, that is 
		$T_2^1$	, $T_2^3$, $T_4^1$ and 	$T_4^3$ are fixed points.\\
		From remark \ref{rem-renumber}, we rewrite $g_1^3$ as $\tilde{g}_3=(1,7,9,15)(3,5,11,13)(2,6,10,14)(4,8,12,16)$ and $g_4^2$ as $\tilde{g}_{14}=(1,4,3,2)(5,16,7,14)(9,12,11,10)(13,8,15,6)$. 
\end{ex}
The same method applies for the computation of  the cycle decomposition of  the bijections $f_i^k$  from the the cycle decompositions of  $\gamma_i$ and $\gamma_k$. We now turn to the study of irretractable and multipermutation solutions.
 \begin{lem}\label{lem-new--irretractable}
	Let $(X,r)$ be a non-degenerate and involutive set-theoretic solution, with $\mid X\mid=n$.  Let  the pair $(\mathcal{X}^2,  \tilde{r})$,  and the functions 
	$	g_i^k\,,\, f_j^l\,: \mathcal{X}^2  \rightarrow\mathcal{X}^2$ be defined as in Definition \ref{defn-new-solution}. Then 
	\begin{enumerate}[(i)]
		\item  If $(X,r)$  is irretractable, then $(\mathcal{X}^2,  \tilde{r})$ is also irretractable.
			\item  If $(X,r)$  is  a multipermutation solution of level $\ell$, then $(\mathcal{X}^2,  \tilde{r})$ is also a multipermutation solution of level $\ell$.
		\item  If $(X,r)$  is  a decomposable solution, then $(\mathcal{X}^2,  \tilde{r})$ is also decomposable.
		\end{enumerate}
	\end{lem}
\begin{proof}
	$(i)$ By contradiction, assume  that for every 
	 $1 \leq i,j,k,l \leq n$,   $g_i^k=g_j^l$. So,    for every $1 \leq m,s \leq n$, $g_i^k(T_s^m)= g_j^l(T_s^m)$. That is, $T_{\sigma_i(s)}^{\sigma_k(m)}=T_{\sigma_j(s)}^{\sigma_l(m)}$, which implies 
	 $\sigma_i(s)=\sigma_j(s)$ and $\sigma_k(m)=\sigma_l(m)$, for every $1 \leq m,s \leq n$. This contradicts  that $(X,r)$  is irretractable.\\
		$(ii)$ We prove that  $\mid\operatorname{Ret}^{l}(\mathcal{X}^2,  \tilde{r})\mid\,=\,\mid\operatorname{Ret}^{l}(X,r)\mid$.  This will imply that if 	$(X,r)$  is  a multipermutation solution of level $\ell$, that is $\mid\operatorname{Ret}^{\ell}(X,r)\mid= 1$, then  $(\mathcal{X}^2,  \tilde{r})$ is  also a multipermutation solution of level $\ell$.	For every $1 \leq i,j \leq n$,  $i \neq j$, whenever  $\sigma_i=\sigma_j$, then for every $1 \leq k\leq n$,  $g_i^k=g_j^k$ and  $g_k^i=g_k^j$. In particular, by substituting $k=i$ or $k=j$, we get  $g_i^i=g_i^k=g_j^k=g_k^i=g_k^j=g_j^j$, that is all the bijections $g_.^*$ with either index $i$ or $j$ are equal.  So, if $x_i \sim x_j$, that is $x_i,x_j \in [x_i]$, then $T_i^i\sim T_i^k\sim T_j^k\sim T_k^i\sim T_k^j\sim T_j^j$, that is these $2n$ elements belong to the same equivalence class that we denote by $[T_i]$. So, there exists a bijection between 	
	$\operatorname{Ret}^{1}(X,r)$ and $\operatorname{Ret}^{1}(\mathcal{X}^2,  \tilde{r})$, and iteratively between 	$\operatorname{Ret}^{l}(X,r)$ and $\operatorname{Ret}^{l}(\mathcal{X}^2,  \tilde{r})$.\\
	$(iii)$ If $(X,r)$  is decomposable, then,  from Proposition \ref{prop-etingof-indecomposable},  there exist $1 \leq s,s'\leq n$, such that for every natural number $m$ and every $i_1,..,i_m$ (not necessarily distinct), $\sigma_{i_1}\sigma_{i_2}...\sigma_{i_m}(s)\neq s'$. This implies that for every $1 \leq k_1,..,k_m\leq n$, $(g_{i_1}^{k_1})\,(g_{i_2}^{k_2})...(g_{i_k}^{k_m})(T_{s}^{.})\neq\, T_{s'}^{*}$. That is, no elements of the form $T_{s}^{.}$ and $T_{s'}^{*}$   belong to the same orbit via the action of  $G(\mathcal{X}^2,  \tilde{r})$ on $\mathcal{X}^2$.  So, $(\mathcal{X}^2,  \tilde{r})$ is   decomposable.
\end{proof}
Note that the proof of Lemma \ref{lem-new--irretractable}$(ii)$ could be used to prove 
Lemma \ref{lem-new--irretractable}$(i)$ also. From Lemma \ref{lem-new--irretractable}$(iii)$, the new solution obtained from  a decomposable solution is  also decomposable. So, a natural question that arises is whether the new solution obtained from  an indecomposable solution is  always  indecomposable. The answer is negative, as Example \ref{ex-continued-n=2} illustrates it. More generally, it is not difficult to show that if  $(X,r)$  is  an indecomposable  permutation solution, then $(\mathcal{X}^2,  \tilde{r})$ is decomposable and furthermore it is  the union of $n$   permutation solutions equivalent to $(X,r)$. So, we  ask when the new solution obtained from  an indecomposable solution is    indecomposable and if there is a criterion on the original solution which permits to decide that. To answer these questions,  we use the following tool. We define  a  table,  $\mathcal{T}$,  for each solution $(X,r)$ in the following way: in the $\ell$-th column we write all the elements $s$  in $X$ for which there exist  $1 \leq i_1,..,i_\ell \leq n$,	such that $\sigma_{i_1}\sigma_{i_2}..\sigma_{i_\ell}(1)=s$. Here are two examples of tables.

	\begin{table*}[h] 
	\begin{tabular}{ |p{1.2cm}|p{1.2cm}|p{1.2cm}|p{1.2cm}|  }
		\hline
			\multicolumn{4}{|c|}{Orbit of $1$  at step $\ell$} \\
			\hline
		 $\ell=1$ & 	 $\ell=2$  & 	 $\ell=3$ &  $\ell=4$ \\
		\hline
	2 &1& 2
		& 1
		\\
		\hline
		\end{tabular}	\hspace{1cm}\begin{tabular}{ |p{1.2cm}|p{1.2cm}|p{1.2cm}|p{1.2cm}|  }
			\hline
			\multicolumn{4}{|c|}{Orbit of $1$  at step  $\ell$} \\
			\hline
		 $\ell=1$ & $\ell=2$  & 	 $\ell=3$ &  $\ell=4$ \\
			\hline
			1,2,3,4 &1,2,3,4 & 1,2,3,4 
			& 1,2,3,4 
			\\
			\hline
		\end{tabular}\\
	
	\caption{At left, the table for the solution in Example \ref{ex-permut-sol-2} and   at right the table for the solution in Example \ref{exemple-irretractable+indecomposable4}.}
\end{table*}

\begin{lem}\label{lem-new-indecomp-iff}
		Let $(X,r)$ be a non-degenerate and involutive set-theoretic solution, with $\mid X\mid=n$. Let $\mathcal{T}$ the table of $(X,r)$.  Let  the pair $(\mathcal{X}^2,  \tilde{r})$,  and the functions 
		$	g_i^k\,,\, f_j^l\,: \mathcal{X}^2  \rightarrow\mathcal{X}^2$ be defined as in Definition \ref{defn-new-solution}.  Let  $\operatorname{Orb}(T_1^1) $ denote  the orbit of $T_1^1$  via the  action of  $G(\mathcal{X}^2,  \tilde{r})$ on $\mathcal{X}^2$.
	 The following statements are equivalent:
	\begin{enumerate}[(i)]
		\item  $(\mathcal{X}^2,  \tilde{r})$  is indecomposable.
		\item $\operatorname{Orb}(T_1^1) \,=\,\mathcal{X}^2$.
		\item For every $1 \leq s,m \leq n$, there exists a natural number $\ell$  and $1 \leq i_1,..,i_\ell, k_1,..,k_\ell \leq n$, such that 
		$\sigma_{i_1}\sigma_{i_2}..\sigma_{i_\ell}(1)=s$ and 
		$\sigma_{k_1}\sigma_{k_2}..\sigma_{k_\ell}(1)=m$. 
		\item For every $1 \leq s,m \leq n$, there exists a natural number $\ell$ such that $s,m$ appear  in column $\ell$ in $\mathcal{T}$.
	\end{enumerate}  
	\end{lem}
\begin{proof}
	From Proposition \ref{prop-etingof-indecomposable}, the solution   $(\mathcal{X}^2,  \tilde{r})$  is indecomposable if and only if the action of  $G(\mathcal{X}^2,  \tilde{r})$ on $\mathcal{X}^2$ is transitive, which is equivalent to  $\operatorname{Orb}(T_1^1) \,=\,\mathcal{X}^2$. This is equivalent to: for every $1 \leq s,m \leq n$, there exists 
	a natural number $\ell$  and $1 \leq i_1,..,i_\ell, k_1,..,k_\ell \leq n$, such that 
	$g_{i_1}^{k_1}g_{i_2}^{k_2}..g_{i_\ell}^{k_\ell}(T_1^1)=T_s^m$,  which is equivalent to 	$\sigma_{i_1}\sigma_{i_2}..\sigma_{i_\ell}(1)=s$ and 
	$\sigma_{k_1}\sigma_{k_2}..\sigma_{k_\ell}(1)=m$,  from Equation (\ref{eqn-defn-new-g}).  From the definition of $\mathcal{T}$, this means that  for every $1 \leq s,m \leq n$, there is a  column $\ell$ in  $\mathcal{T}$ in which both $s$ and $m$ appear together.
\end{proof}
We are  now  able to state condition  	$( \mathfrak{C})$ in two equivalent ways, to formulate Theorem \ref{thmA-construction} in a more precise way and to give its proof.
\begin{thm}\label{thm-construction}
	Let $(X,r)$ be a non-degenerate and involutive set-theoretic solution, with $\mid X\mid=n$ and table $\mathcal{T}$. Let  the pair $(\mathcal{X}^2,  \tilde{r})$,  and the functions 
	$	g_i^k\,,\, f_j^l\,: \mathcal{X}^2  \rightarrow\mathcal{X}^2$ be defined as in Definition \ref{defn-new-solution}.  Then, for every natural number $k \geq 1$,    there exists a  non-degenerate and involutive set-theoretic solution, $(\mathcal{X}^2,  \tilde{r})$,  of size $n^{2^k}$ induced from $(X,r)$. Furthermore, the following hold:
	\begin{enumerate}[(i)]
		\item If  $(X,r)$  is irretractable or a multipermutation solution of level $\ell$, then 
		$(\mathcal{X}^2,  \tilde{r})$ is  irretractable or  a multipermutation solution of level $\ell$, respectively.
		\item 	If  $(X,r)$  is decomposable, then 
		$(\mathcal{X}^2,  \tilde{r})$ is also decomposable.
			\item 	If  $(X,r)$  is of class $m$, then 
		$(\mathcal{X},  \tilde{r})$ is also  of class $m$.
	\end{enumerate} 
Assume that $(X,r)$  is indecomposable and satisfies the following  condition 
	$( \mathfrak{C})$ written in two equivalent forms:\\
	
	$( \mathfrak{C})$  There exists 
	a natural number $\ell$ such that all the elements of $X$ appear in column $\ell$. \\
	$( \mathfrak{C})$  There exists 
	a natural number $\ell$ such that 
	for every $1 \leq s \leq n$,  there exist $1 \leq i_1,..,i_{\ell}\leq n$   such that 
	$\sigma_{i_1}\sigma_{i_2}..\sigma_{i_{\ell}}(1)=s$. \\ 
	 Then   $(\mathcal{X}^2,  \tilde{r})$  is also  indecomposable and satisfies   condition  $( \mathfrak{C})$.
	Moreover, the same process can be repeated iteratively to obtain  indecomposable solutions of size $n^{2^k}$,  for every natural number  $k \geq 1$.
\end{thm}
\begin{proof}
	From Lemma \ref{lem-new--invol+braided},  $(\mathcal{X}^2,  \tilde{r})$ 
	is a non-degenerate and involutive set-theoretic solution, with $\mid \mathcal{X}^2\mid=n^2$.  $(i)$ and $(ii)$ result from Lemma \ref{lem-new--irretractable} and $(iii)$ from Lemma \ref{lem-new-class}.  If $(X,r)$  is indecomposable and satisfies 
	$( \mathfrak{C})$, then  from Lemma \ref{lem-new-indecomp-iff}$(iv)$, $(\mathcal{X}^2,  \tilde{r})$  is  indecomposable. Indeed,  condition 
	$( \mathfrak{C})$  implies condition $(iv)$ from  Lemma \ref{lem-new-indecomp-iff}.  We show that $(\mathcal{X}^2,  \tilde{r})$   satisfies condition  $( \mathfrak{C})$, that is   there exists 
	a natural number $\ell'$ such that 
	for every $1 \leq s,m \leq n$,  there exist $1 \leq i_1,..,i_{\ell'}, k_1,..,k_{\ell'} \leq n$   such that 
	$g_{i_1}^{k_1}g_{i_2}^{k_2}..g_{i_{\ell'}}^{k_{\ell'}}(T_1^1)=T_s^m$,   which is equivalent to 	$\sigma_{i_1}\sigma_{i_2}..\sigma_{i_{\ell'}}(1)=s$ and 
	$\sigma_{k_1}\sigma_{k_2}..\sigma_{k_{\ell'}}(1)=m$.  This last statement is true, with $\ell'=\ell$,  since  $(X,r)$ satisfies condition 	$( \mathfrak{C})$. Clearly,  iteratively doing  the same process, we  can  obtain  indecomposable solutions of size $n^{2^k}$, for every natural number  $k \geq 1$.
\end{proof}
Clearly,  condition 
$( \mathfrak{C})$  implies condition $(iv)$ from  Lemma \ref{lem-new-indecomp-iff}. We do not know if in our context these two conditions are equivalent. 
\section{Applications to cryptography}
In this paper, based on the  construction from Section \ref{subsec-iterative}, we suggest a public key encryption method and a signature  method. Using the same ideas, it might be possible to construct a method for the authentication and the key exchange procedures also. We first present our methods, then   evaluate their complexity and at last examine their strengths and weaknesses.

\subsection{Suggestion of a public key encryption method and a signature method}
Let  $(X,r)$ be a  non-degenerate and involutive set-theoretic solution of the   YBE, with $\mid X\mid=n$.  From Theorem \ref{thmA-construction},   for every $k\geq 1$, there exists a  non-degenerate and involutive set-theoretic solution of size $n^{2^k}$. We call \emph{pumping up a solution} the process of  constructing a solution of size $n^{2^k}$  from a solution  satisfying the above  conditions.  Such   a solution  is called \emph{a pumped-up solution}, and we denote by $\hat{g}_1,...,\hat{g}_{n^{2^k}}$ the permutations (of kind $\sigma$) defining  it; $(X,r)$   is called \emph{the original solution}. Note that the original solution might be a solution which  is itself the result of some  previous pumping up.
 From Corollary \ref{cor1}, if  $(X,r)$  is 	irretractable, indecomposable and satisfies   the  condition $( \mathfrak{C})$,  then all  the pumped-up solutions are also  irretractable and  indecomposable. \\
 \textbf{Our suggestion of a public key encryption method:}\\
 
 Let $M$ be the message Bob wants to send to Alice. 
 \begin{itemize}
 	\item Begin with an original solution and pump it to get a solution of size $n^{2^k}$, for  $k$ large enough. 
 	\item  Make $2^k$ (or $k$) public. Although  $2^k$  is  public, the original solution, its size, the new solution and its  size  are  effectively hidden from everyone else.
 	\item   Choose $ 1 \leq i  \leq n^{2^k}$ such that  $i << n^{2^k}$ and make it public.
 	\item Bob computes the permutation $\hat{g}_i$ in $S_{n^{2^k}}$, the symmetric group on  $n^{2^k}$ elements, from the original solution.
 	\item  Bob breaks the message $M$ into a series of blocks such that it is possible to  represent each block as  an integer which is smaller or equal to $n^{2^k}$. For that, it is possible to  use any standard representation. The purpose here is not to encrypt the message but only to get it into the numeric form necessary for encryption.
  \item Bob encrypts  the message by applying the permutation  $\hat{g}_i$ on each numerical block. That is, the result (the ciphertext $C$) is the image of each numerical block by $\hat{g}_i$. 
  \item Alice computes  the permutation  $\hat{g}^{-1}_i$ in $S_{n^{2^k}}$, from the original solution.
  \item Alice decrypts the ciphertext $C$ by applying the permutation  $\hat{g}^{-1}_i$ on each numerical block. 
 \end{itemize}

Let illustrate with a small  example:  the original solution is from Example \ref{exemple-irretractable+indecomposable4} and it is pumped up    to a solution of size $n^{2^2}$, where $n=4$.  A part of the computations  is done in Example \ref{ex-pump4}. 
\begin{ex}\label{ex-wonderful}
		 The public key is  $2^2$ and make  $i=46$  public. Bob wants to send to Alice the title of a song of Black (late 90'), so his  message  is:\\
	ITS A WONDERFUL LIFE\\
He  can encode two letters per block, substituting a two-digit number for each letter:\\ blank $= \,00$, $A = \,01$, $B =\, 02$, ..., $Z =\, 26$.  In order to have each numerical block less or equal 256 and of the same length,  the message is encoded as:\\
09 20 19 00 01 00 
 23 15 14 04 05 18 06 21 12 00 12 09 06 05\\
Bob needs to compute the permutation $\hat{g}_{46}$ in $S_{256}$ or at least  the result of the application of this permutation on each of the numbers above. We will discuss later on how to compute effectively $\hat{g}_{46}$ in $S_{256}$.
We extend the domain of definition of $\hat{g}_{46}$, so that $\hat{g}_{46}(00)=00$, and  
the  whole message is  then enciphered as:\\
108  83 82 00 100 00 94 102 101 99  112 81 109 96  107 00 107 108 109 112\\
In order to decipher the message, Alice computes $\hat{g}^{-1}_{46}$ in $S_{256}$ or at least  the result of the application of this permutation on each of the numbers in the enciphered message and retrieves the original message.
Maybe it is possible to only require   each numerical block less or equal 256 and remove the requirement that each numerical block has  the same length. 
\end{ex}

\begin{rem}
	\begin{enumerate}
		\item  Applying first  $\hat{g}^{-1}_i$ on the numerical blocks representing the message $M$ and next $\hat{g}_i$ on the result produces back the numerical blocks representing $M$. That is, if a message $M$ is first deciphered and then enciphered, $M$ is the result. 
		\item Given the public key $i$, Alice knows the  encryption procedure of Bob (the permutation $\hat{g}_{i}\in S_{n^{2^k}}$) and Bob can  also compute the decryption procedure of Alice ($\hat{g}_{i}^{-1}\in S_{n^{2^k}}$).
	\end{enumerate}
\end{rem}

 \textbf{Our suggestion of a signature method:}\\

Let $M$ denote  the message Bob wants to send to Alice and $S$ his signature.
\begin{itemize}
	\item Begin with an original solution and pump it to get a solution of size $n^{2^k}$, for  $k$ large enough. 
	\item  Make $2^k$ (or $k$) public. Although  $2^k$  is  public, the original solution, its size, the new solution and its  size  are  effectively hidden from everyone else.
	\item   Alice chooses $ 1 \leq i  \leq n^{2^k}$ such that  $i << n^{2^k}$ and makes it public.
	\item   Bob chooses $ 1 \leq j \leq n^{2^k}$ such that  $j << n^{2^k}$ and makes it public.
	\item Alice  computes the permutations $\hat{g}_j$  and $\hat{g}^{-1}_i$ in $S_{n^{2^k}}$,  from the original solution.
	\item  Bob computes the permutations $\hat{g}_i$  and $\hat{g}^{-1}_j$ in $S_{n^{2^k}}$, from the original solution.
	\item  Bob breaks the message $M$ into a series of blocks such that it is possible to  represent each block as  an integer which is smaller or equal to $n^{2^k}$. For that, it is possible to  use any standard representation. The purpose here is not to encrypt the message but only to get it into the numeric form necessary for encryption.
	\item Bob computes $S=\,\hat{g}^{-1}_j(M)$,  then encrypts $S$ by applying the permutation  $\hat{g}_i$ on $S$ and sends $\hat{g}_i(S)\,=\,\hat{g}_i\hat{g}^{-1}_j(M)$ to Alice. 
	\item Alice computes  $S$ by applying  the permutation  $\hat{g}^{-1}_i$ on $\hat{g}_i(S)$. 
	\item Alice retrieves  $M$ by applying the permutation  $\hat{g}_j$ on $S$:
	$\hat{g}_j(S)\,=\,\hat{g}_j\hat{g}^{-1}_j(M)\,=\,M$. 
\end{itemize}
In both methods, the reason for the choice of  the public key  $ 1 \leq i  \leq n^{2^k}$ such that  $i << n^{2^k}$ is to avoid giving an intruder any clue on $ n^{2^k}$,  the size of the chosen pumped-up solution.  Let illustrate with a small  example:  the original solution is from Example \ref{exemple-irretractable+indecomposable4} and it is pumped up    to a solution of size $n^{2^2}$, where $n=4$.  A part of the computations  is done in Examples  \ref{ex-pump4} and \ref{ex-wonderful}. 
\begin{ex}\label{ex-wonderful2}
	The public key is  $2^2$.  Alice chooses $i=46$  and makes it    public. Bob chooses $j=3$  and makes it    public. Bob wants to send to Alice  a signed message,  his  message  $M$ is the title of a song of Black (late 90'):\\
	ITS A WONDERFUL LIFE\\
The  message is encoded as:\\
	09 20 19 00 01 00 
	23 15 14 04 05 18 06 21 12 00 12 09 06 05\\
	Bob needs to compute  the signature $S$: he computes the permutation  $\hat{g}_{3}^{-1}$ in $S_{256}$ or at least  the result of the application of this permutation on each of the numbers above. He applies  $\hat{g}_{3}^{-1}$ on the message and obtains $S=\hat{g}_{3}^{-1}(M)$:\\
247 208 205  00 255 00 193 249 250 256 243 206 242 195 248  00 248 247 242 243 \\
He then encrypts his signature: he applies  $\hat{g}_{46}$  on $S$ and sends to Alice:\\
62 39 40 00 54 00 36 60 57 55 50 37 49 34 63 00 63 62 49 50\\
 To obtain $S$, Alice  applies  $\hat{g}^{-1}_{46}$  on each of these numbers. She  then retrieves the original message $M$ by applying $\hat{g}_{3}$ on $S$.
\end{ex}
Note that in the suggested method, given $i$, Bob knows both the encryption and decryption procedures of Alice ($\hat{g}_i$  and $\hat{g}^{-1}_i$ respectively) and given $j$,  Alice  knows both the encryption and decryption procedures of Bob  ($\hat{g}_j$  and $\hat{g}^{-1}_j$). This is not the case in standard signature procedures.
 \subsection{Evaluation  of the complexity of the suggested methods} \label{sec_tree}

The first question that arises about both  methods is the question of their   security. However, before we get into the examination of this question, we need to know if the encryption and decryption procedures in the first method and the keys computations in the second method can be done  sufficiently easily.   Both methods rely  on the computation of one or several  permutations  $\hat{g}_i$ in  the symmetric group on  $n^{2^k}$ elements, based on the knowledge of the original solution. So, we need to understand how to do this calculation and what is its complexity. We present a way to compute $\hat{g}_i$. Yet,  there might be  more efficient ways to do it.

As a first step, we  construct a finite regular binary tree $\mathscr{T}_i$ in the following way:\\
\begin{gather*}
\textrm{The  top level is level} \;k \,\textrm{and the root  at  level} \;k \;\textrm{is labelled  by} \;i.\\
\textrm{At level}\; k-1:  \textrm{the leftmost node is} \; \alpha_0=\,\lceil \frac{i}{n^{2^{k-1}}}\rceil\;
\textrm{and the rightmost node is} \;\alpha_1\equiv i \, (mod\, n^{2^{k-1}})\\
 \textrm{At level} \; k-2: \textrm{from left}\; \alpha_{0,0}=\,\lceil \frac{\alpha_0}{n^{2^{k-2}}}\rceil,\;
 \alpha_{0,1}\equiv \,\alpha_0(mod \,n^{2^{k-2}}),\;
 \alpha_{1,0}=\,\lceil \frac{\alpha_1}{n^{2^{k-2}}}\rceil, \;
 \textrm{} \;
  \alpha_{1,1}\equiv \,\alpha_1(mod\, n^{2^{k-2}})\\
  \textrm{At level} \; k-m: \textrm{inductively}\;
   \alpha_{\operatorname{ad}_{\{k-m+1\},\,0}}=\,\lceil
   	 \frac{\alpha_{\operatorname{ad}_{\{k-m+1\}}}}{n^{2^{k-m}}}\rceil\;
   	  \textrm{and} \;	\alpha_{\operatorname{ad}_{\{k-m+1\}},\,1}\equiv \,\alpha_{\operatorname{ad}_{\{k-m+1\}}}(mod \,n^{2^{k-m}})\\
\textrm{where}\;   \alpha_{\operatorname{ad}_{\{k-m+1\}}}\; \textrm{denotes the label at some  address at level}\; k-m+1, \; \textrm{and} \;m\leq k.
      \end{gather*}
      At  level $0$,  the bottom level,  there are $2^k$ nodes and each pair of sibling nodes  is connected to a unique parent at level 1. So,  each pair of sibling nodes at level $0$ can be described by  $(\alpha_{\operatorname{ad}_{\{1\}},\,\tiny 0 }\,,\, \alpha_{\operatorname{ad}_{\{1\}},\,1})$,  where $\operatorname{ad}_{\{1\}}$ denotes  their parent's address from  level $1$. Inductively,  $\alpha_{\operatorname{ad}_{\{1\}}}$ can be described by $\alpha_{\operatorname{ad}_{\{2\}},\,\tiny 0 }$ or  $\alpha_{\operatorname{ad}_{\{2\}},\,1}$,  where $\operatorname{ad}_{\{2\}}$ denotes its  parent's address from  level $2$, and so on until level $k$. We illustrate the construction of the tree  for $i=46$ from Example \ref{ex-wonderful} ($n=4$, $k=2$, $n^{2^k}=256$, $n^{2^{k-1}}=16$).  
	\begin{figure}[H]\label{fig-tree}
  			\begin{tikzpicture}[level distance=1.5cm,
  			level 1/.style={sibling distance=6.4cm},
  			level 2/.style={sibling distance=4cm},
  			level 3/.style={sibling distance=4cm}]
  			\node {$\bf{i=46}$}
  			child {node {$\lceil \frac{46}{4^{2^{k-1}}}\rceil=\bf{3}$}
  			child {node {$\lceil \frac{3}{4^{2^{k-2}}}\rceil=\bf{1}$}}  child {node {$\bf{3} \,$$\equiv 3$}}
  				}
  				child {node {$\bf{14}\,$$\equiv 46\mod 4^{2^{k-1}}$}
  				child {node {$\lceil \frac{14}{4^{2^{k-2}}}\rceil=\bf{4}$} 
  				}
  				child {node {$\bf{2}\,$$ \equiv 14\mod 4^{2^{k-2}}$}}
  			};
  			\end{tikzpicture}
  		\caption{The binary tree $\mathscr{T}_{46}$ for $i=46$,  and $n^{2^k}$, with $n=4$ and $k=2$.}
  	\end{figure}
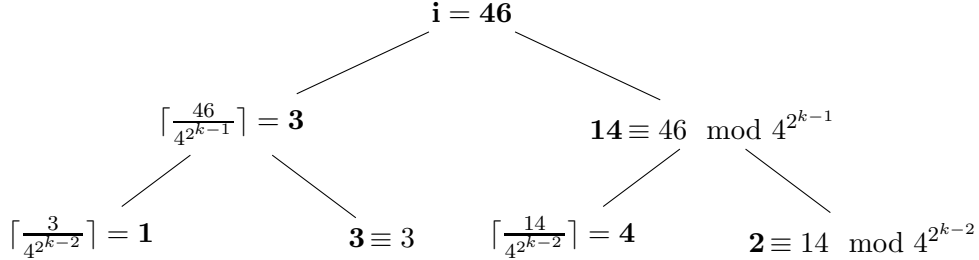

After the construction of $\mathscr{T}_i$, we compute $\hat{g}_i$ from it  in the following way:\\
We begin at level $0$, the bottom level, and for each pair of sibling nodes with labels\\ $(\alpha_{\operatorname{ad}_{\{1\}},\,{\tiny 0 }  }\,,\, \alpha_{\operatorname{ad}_{\{1\}},\,1})$,  where $\operatorname{ad}_{\{1\}}$ denotes their parent's address from level $1$, we compute the following permutation in $S_{n^{2}}$:
 \begin{equation*}
\tilde{g}_{\alpha_{\operatorname{ad}_{\{1\}}}}\,\in S_{n^{2}}\,:\;\;\;
\tilde{g}_{\alpha_{\operatorname{ad}_{\{1\}}}}\,=\,g_{j}^{l},  \;\textrm{where} \;\;j=\alpha_{\operatorname{ad}_{\{1\}},\,0}\;\textrm{ and }\;l=\alpha_{\operatorname{ad}_{\{1\}},\,1}
\end{equation*}
The computation is done using the permutations $\sigma_j$ and $\sigma_l$  from the original solution. Each computed permutation  $\tilde{g}_{\alpha_{\operatorname{ad}_{\{1\}}}}$ is attached to the corresponding node at level $1$.  At this level, there are $2^{k-1}$ computed permutations in $S_{n^{2}}$.\\
Next, we move to level $1$,  
and for each pair of sibling nodes with labels $(\alpha_{\operatorname{ad}_{\{2\}},\,\tiny 0 }\,,\, \alpha_{\operatorname{ad}_{\{2\}},\,1})$,  where $\operatorname{ad}_{\{2\}}$ denotes their parent's address from level $2$, we compute the following permutation in $S_{n^{2^2}}$:
\begin{equation*}
\tilde{\tilde{g}}_{\alpha_{\operatorname{ad}_{\{2\}}}}\,\in S_{n^{2^2}}\,:\;\;\tilde{\tilde{g}}_{\alpha_{\operatorname{ad}_{\{2\}}}}\,=\,\tilde{g}_{s}^{m} \;\textrm{where} \;\;s=\alpha_{\operatorname{ad}_{\{2\}},\,0}\;
\textrm{and}\;m=\alpha_{\operatorname{ad}_{\{2\}},\,1}
\end{equation*} 
The computation is done using the permutations $\tilde{g}_{\alpha_{\operatorname{ad}_{\{2\}},0}}$ and $\tilde{g}_{\alpha_{\operatorname{ad}_{\{2\}},1}}$  from $S_{n^{2}}$, computed at the previous step. 
Each computed permutation  $\tilde{\tilde{g}}_{\alpha_{\operatorname{ad}_{\{2\}}}}$ is attached to the corresponding node at level $2$. At this level, there are $2^{k-2}$ computed permutations in $S_{n^{2^2}}$. Next, we move to level $2$,  and repeat the same process until we arrive to level $k-1$ and obtain $\hat{g}_i=\tilde{\tilde{g}}_{\alpha_{0}}^{\alpha_{1}}$.  The total number of  permutations computations  is $2^{k-1}+2^{k-2}+...+2+1$, that is $2^{k}-1$, from the formula of the sum of a geometric sequence. 

\begin{ex}
We illustrate the computation of $\hat{g}_{46}$ using  $\mathscr{T}_{46}$ from Figure \ref{fig-tree}. We begin at level $0$. We compute $\tilde{g}_{3}=g_1^3$, since $\alpha_{0}=3$,  $\alpha_{0,0}=1$,  $\alpha_{0,1}=3$,  
and compute $\tilde{g}_{14}=g_4^2$, since $\alpha_{1}=14$, $\alpha_{1,0}=4$
and $\alpha_{1,1}=2$. From Example \ref{ex-pump4}, \\ $\tilde{g}_3=(1,7,9,15)(3,5,11,13)(2,6,10,14)(4,8,12,16)$ and \\ $\tilde{g}_{14}=(1,4,3,2)(5,16,7,14)(9,12,11,10)(13,8,15,6)$. \\
We move to level $1$ (which is also level $k-1$ here), and compute $\hat{g}_{46}$ in $S_{256}$, by  $\hat{g}_{46}=\tilde{g}_{3}^{14}$, since 
 $\alpha_{0}=3$,  $\alpha_{1}=14$.  As an example, $\hat{g}_{46}(09)=108$:  $09$ corresponds to an element of the form $T_1^9$,  which image is $T_{7}^{12}$, since $\tilde{g}_3(1)=7$ and $\tilde{g}_{14}(9)=12$, renumbering it gives $108=6\cdot16+12$. Another example $\hat{g}_{46}(23)=94$:  $23$ corresponds to an element of the form $T_2^7$,  which image is $T_{6}^{14}$, since $\tilde{g}_3(2)=6$ and $\tilde{g}_{14}(7)=14$, renumbering it gives $94=5\cdot16+14$. 
\end{ex}
From the above,  given $i$, the computation of $\hat{g}_i$ in $S_{n^{2^k}}$ requires:\\
\begin{itemize}
	\item $2^k-1$ divisions.
	\item $2^k-1$ residues computations.
	\item $2^k-1$ permutations computations.
	\item enough space to keep  a binary tree with $2^{k+1}-1$ labelled nodes. 
\end{itemize}
The computation of a  permutation  in   $S_{n^{2^l}} $ requires  $n^{2^{l}}$ substitutions, and $n^{2^l}$ simple operations to renumber each element of the form $T_i^k$, as described in Remark \ref{rem-renumber}.  So,  for every $1 \leq l \leq k$, the total number of operations is   $2 \cdot((2^{k-1})\cdot n^{2^1}+(2^{k-2})\cdot n^{2^2}+...+(2^{k-k})\cdot n^{2^k})$, that is   $2^{k}\cdot n^{2^1}+2^{k-1}\cdot n^{2^2}+...+2\cdot n^{2^k}$ operations.  

We turn now to the evaluation  of the complexity of the suggested methods. For our suggestion of a public key encryption method, both Alice and Bob need for one value $i<<n^{2^k}$ to compute $\hat{g}_i$ in $S_{n^{2^k}}$, or at least  the result of the application of this permutation on each of the needed numbers. From the list above,  assuming that it takes $10^{-9}$ seconds for a computer to  make an operation, as a gross approximation:
 \begin{gather}
\nonumber\textrm{Computing} \;\hat{g}_i \;\textrm{in} \;S_{n^{2^k}}\; \textrm{requires}\;\;\approx\\
((2^{k+1}-2)\,+\,(2^{k}\cdot n^{2^1}+2^{k-1}\cdot n^{2^2}+...+2\cdot n^{2^k}))\cdot10^{-9} \;\textrm{seconds}
 \end{gather}\label{eqn-time-pump1}
Note that to compute  the result of the application of  $\hat{g}_i$ on a single number or some numbers  spares us only one permutation computation, the computation of $\hat{g}_i$ itself, and  all the previous permutations are needed. Yet, this is an important "economy", if $k$ is large, since  the computation of $\hat{g}_i$ itself  requires  $2\cdot n^{2^k}$ computations. So, a possible option is to stop  the computations  one step before  the computation of  $\hat{g}_i$ itself, and at each message compute only  the result of the application of  $\hat{g}_i$  on each of the needed numbers.

	\begin{table}[h] 
	\centering
	\begin{tabular}{ |p{2cm}|p{2cm}|p{2cm}|p{2cm}|p{2cm}| }
		\hline
		&$k=2$ &	$k=3$ &$ k=4$ &$k=5$ \\
		\hline
		$n=2$ &\cellcolor{yellow}$5\cdot10^{-8}$&\cellcolor{yellow} $6\cdot10^{-7}$ &\cellcolor{yellow}$1\cdot10^{-4}$&\cellcolor{yellow}$8.6$
		\\
		\hline
		$n=3$ &\cellcolor{yellow}$2\cdot10^{-7}$&\cellcolor{yellow} $1\cdot10^{-5}$ &\cellcolor{yellow}$0.08$&$3\cdot10^{6}$
		\\
		\hline
		$n=4$&\cellcolor{yellow}$5\cdot10^{-7}$& \cellcolor{yellow}$1\cdot10^{-4}$ &\cellcolor{yellow}$8.5$&$3\cdot10^{10}$\\
		\hline
		$n=5$ &\cellcolor{yellow}$1\cdot10^{-6}$&\cellcolor{yellow}$7\cdot10^{-4}$  & \cellcolor{yellow}$305$ &$...$
		\\
			\hline
		$n=6$ &\cellcolor{yellow}$2\cdot10^{-6}$&\cellcolor{yellow}$0.003$  & \cellcolor{yellow}$5642$ &$...$
		\\
			\hline
		$n=7$ &\cellcolor{yellow}$5\cdot10^{-6}$&\cellcolor{yellow}$0.01$  & \cellcolor{yellow}$6\cdot10^{4}$ &$...$
		\\
			\hline
		$n=8$ &\cellcolor{yellow}$8\cdot10^{-6}$&\cellcolor{yellow}$0.03$  & $5\cdot10^{5}$ &$...$
		\\
			\hline
		$n=9$ &\cellcolor{yellow}$1\cdot10^{-5}$&\cellcolor{yellow}$0.08$  & $3\cdot10^{6}$ &$...$
		\\
			\hline
		$n=10$ &\cellcolor{yellow}$2\cdot10^{-5}$&\cellcolor{yellow}$0.2$  & $2\cdot10^{7}$ &$...$
		\\
			\hline
		$n=16$ &\cellcolor{yellow}$1\cdot10^{-4}$&\cellcolor{yellow}$8.5$  & $3\cdot10^{10}$ &$...$
		\\
			\hline
		$n=25$ &$\cellcolor{yellow}7\cdot10^{-4}$&\cellcolor{yellow}$305$  & $...$ &$...$
		\\
			\hline
		$n=36$ &\cellcolor{yellow}$0.003$&\cellcolor{yellow}$5\cdot10^{3}$  & $...$ &$...$
		\\
	\hline
$n=49$ &\cellcolor{yellow}$0.01$&\cellcolor{yellow}$6\cdot10^{4}$  & $...$ &$...$
\\
	\hline
$n=64$ &\cellcolor{yellow}$0.03$&$5\cdot10^{5}$  &$...$ &$...$
\\
	\hline
$n=81$ &\cellcolor{yellow}$0.08$&$3\cdot10^{6}$  & $...$ &$...$
\\	\hline
$n=100$ &\cellcolor{yellow}$0.2$&$2\cdot10^{7}$  & $...$ &$...$
\\
		\hline
	\end{tabular}\\
	\caption{The approximate time in seconds to compute one $\hat{g}_i$ in the pumped-up solution of size $n^{2^k}$.   One year is $\approx 3.1 \cdot10^7$ seconds,  one month  is $\approx 2.6 \cdot10^6$ seconds and  one day is $\approx 8.6 \cdot10^4$ seconds.}\label{table_time_one_gi}
\end{table}
 For our suggestion of a signature method, both Alice and Bob need to compute either two permutations  in $S_{n^{2^k}}$,  or  the result of the application of  these two   permutations on some numbers. 
 From the  above,  assuming that it takes $10^{-9}$ seconds  to  make an operation, this requires  two times the time  in Equation  \ref{eqn-time-pump1}, that is 
$ \approx \,2\cdot \,((2^{k+1}-2)\,+\,(2^{k}\cdot n^{2^1}+2^{k-1}\cdot n^{2^2}+...+2\cdot n^{2^k}))\cdot10^{-9} \;\textrm{seconds}$.\\
Note that in both methods, the computation of  the permutation $\hat{g}_i$ in the pumped-up solution of size $n^{2^k}$ is independent of the message and can be done ahead of time if needed.\\

\begin{figure}[H]\label{fig-tree-ismall}
	\begin{tikzpicture}[level distance=1cm,
		level 1/.style={sibling distance=4cm},
		level 2/.style={sibling distance=2cm},
		level 3/.style={sibling distance=1.5cm}]
		\node {$i=3$}
		child {node {$1$}
			child { node {$1$} child {node {$1$}} child {node {$1$}}
			}
			child { node {$1$} child {node {$1$}} child {node {$1$}}
		}
		}
		child {node {$3$}
			child {node {$1$} child {node {$1$}} child {node {$1$}}
			}
			child { node {$3$} child {node {$1$}} child {node {$3$}}
		}
		};
		\end{tikzpicture}
	\caption{The binary tree  $\mathscr{T}_{3}$ for small $i=3$,   with $n\geq 3$,  $k=3$.} 
\end{figure}
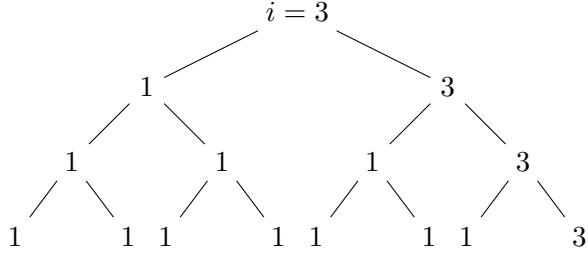

As Table \ref{table_time_one_gi} illustrates it,  even with the knowledge of the original solution,  the time required to compute  one $\hat{g}_i$ in the pumped-up solution of size $n^{2^k}$  grows very fast as $n$ and $k$ grow. For large  $k$, an easy step can be done to reduce the  number  of computations. Indeed, by choosing $i$ very small, that is $i <n$,  one can spare all the computations in the  tree construction and reduce drastically the number of  permutations  computations. Indeed, as Figure  \ref{fig-tree-ismall} illustrates it,  in the binary tree  $\mathscr{T}_{i}$ for  a small $i$, all the nodes are labelled $1$ except the nodes in  the rightmost branch of  $\mathscr{T}_{i}$ which are labelled $i$. So, $2^{k+1}-2$ computations are spared  in the  tree construction. Furthermore, at each level $l<k-1$, instead of $2^{k-1-l} $ permutations computations, there are only $2$ permutations computations. So,  we need $2^{2}\cdot n^{2^1}+2^{2}\cdot n^{2^2}+...+2\cdot n^{2^k}$ computations. 
 \begin{gather}
\nonumber\textrm{Computing} \;\hat{g}_i ,\;\textrm{for small}\; i,\; \textrm{in}\; S_{n^{2^k}}\; \textrm{requires}\;\;\approx\\
2(2 n^{2^1}+2 n^{2^2}+...+n^{2^k})\cdot10^{-9} \;\textrm{seconds}
\end{gather}\label{eqn-time-pump-small-i}

\subsection{Estimation of the security of the suggested methods }

 An intruder who wants to discover the message Bob sends to Alice  has to guess a permutation  from the IYB group in  $S_{n^{2^k}}$,  the symmetric group on  $n^{2^k}$ elements (with $n^{2^k}$ large enough). As  he has no information on the original solution, nor does he know $n$, he needs to search the permutation in $S_{n^{2^k}}$,   for several $n$. So, in the case of a brute force attack, our  intruder has to test every permutation in  the symmetric group on  $n^{2^k}$ elements for every $n \geq 2$ until he breaks it.  We recall that the symmetric group on  $n^{2^k}$ elements contains  $(n^{2^k})!$ permutations. So,   assuming   that it takes $10^{-8}$ seconds for a computer to  make a search, the amount of time required for a computer to  search for all the permutations in $S_{n^{2^k}}$, for a given $n$,  is   $(n^{2^k})!\,\cdot10^{-8}$ seconds.   As an example, if $n=4$ and $k=2$, the order of the symmetric group $S_{256}$ is $256!$, which is approximately $8.1\cdot10^{506}$, using Stirling approximation $m!  \approx \sqrt{2\pi m}\,(\frac{m}{e})^m$. So to search for all the permutations in $S_{256}$ requires $8.1\cdot10^{498}$ seconds, that is many years.\\
 
 Another kind of "smarter"  brute force attack is if our intruder anticipates the possible cycle decompositions of the permutations in all the possible IYB groups of the pumped-up solutions. Indeed, in  the IYB group of the original solution, not all the cycle decompositions  can occur and  the cycle decompositions of the permutations occurring   depend on $n$ and on the class of the solution $m$ (see Remark  \ref{remark_defn_class}$(ii)$).  The cycle decompositions of the permutations in the  IYB group of a  pumped-up solution are determined by the 
cycle decompositions of the permutations in the  IYB groups of the original solution and the intermediate  pumped-up solutions. So, our intruder can instead of searching all the  permutations 
 in $S_{n^{2^k}}$,   for every $n \geq 2$ until he breaks it, can search for several possible    cycle decompositions in each group of permutations.\\
 
 The number of permutations of a given cycle decomposition in $S_k$ is:
 \begin{equation}
 \frac{k!}{\prod\limits_{d=1}^{d=k}(n_d)!\,d^{n_d}}
 \end{equation}\label{eqn-number-permut}
 where $n_d$ is the number of cycles of length $d$.\\
  As an example, if $k=2$, our intruder anticipates that for $n=4$ and certain kind of solutions with $n=4$ (like the solution in Example \ref{ex-pump4}), the  possible permutations in $S_{256}$ have the possible form:  either the disjoint product of  $64$  ($4$-cycles) or the  the disjoint product of  $n_2$  ($2$-cycles) and 
   $n_1$  ($1$-cycles), where $n_1$ and $n_2$ may vary. Then, from Equation \ref{eqn-number-permut}, the number of  permutations in $S_{256}$ which are the  disjoint product of  $64$  ($4$-cycles) is, using Stirling approximation,    approximately $\frac{8\cdot10^{506}}{64! 4^{64}}\approx 2\cdot10^{379}$.  So, to search for all the permutations with this cycle decomposition in $S_{256}$ requires $2\cdot10^{371}$ seconds, that is many years.\\

   Now, assume that  instead of a brute force attack, our intruder decides  to find the pumped-up  solution of size $n^{2^k}$, using the construction we described in the paper.  As he  does not know $n$, he needs to go through all  the possible values of $n$ until he breaks it. Furthermore, for each $n$, he has to apply the construction on  each solution of size $n$, since he has no information on the original solution, that is for each $n$ and for each solution he needs  to compute $\hat{g}_i$ and test it.  For each $n \geq 2$, there are several candidates for the original solution, and this number grows very fast.  The following table describes  the growth of the number of solutions of a given size (the number of solutions up to isomorphism is   from \cite{arxiv-ybe}, and a part of it from \cite{etingof}):
 	\begin{table}[h]
 		 \centering
 	\begin{tabular}{ |p{2cm}|p{0.5cm}|p{0.5cm}|p{0.6cm}|p{0.8cm}|p{1.3cm}| p{1.3cm}| p{1.3cm}| p{1.5cm}| p{1.5cm}|  }
 		\hline
 		$n=$&$2$ &	$3$ & 4 &5&6& 7 &8&9&10\\
 		\hline
 		all	sol. & 2&  12 &168&2640& $8.2 \cdot 10^{4}$ &$2.6 \cdot 10^6$&$1.6\cdot 10^{8}$&-&-\\
 		\hline
 		sol. up iso & 2&  5 &23&88& 595 &3456&34.530&321.931 &4.895.272
 		\\
 		\hline
 	\end{tabular}\\
 		\caption{The approximate total  number of non-degenerate and involutive solutions and  their number up to isomorphism,    for  $n \leq 10$.}\label{table_nb-solutions}
 \end{table}\\
From Table \ref{table_nb-solutions},  the set of solutions from which one can choose the original solution  grows very fast as $n$ grows. In the classification of the non-degenerate and involutive solutions,  the interest is on the number of solutions up to isomorphism. However, in our context, the interest is on the total number of solutions of a given size. Using the data from Table \ref{table_nb-solutions}, we compute the estimated  time required for an intruder to compute $\hat{g}_i$ for all the solutions of size $n$. For $n=9,10$, we assume that the total number of solutions is 1000 times the number of solutions up to isomorphism (which is a lower estimation than the real number).
\begin{table}[h] 
	\centering
	\begin{tabular}{ |p{2cm}|p{2cm}|p{2cm}|p{2cm}|p{2cm}| }
		\hline
		&$k=2$ &	$k=3$ &$ k=4$ &$k=5$ \\
		\hline
		$n=2$ &$1\cdot10^{-7}$& $1\cdot10^{-6}$ &$2\cdot10^{-4}$&$17$
		\\
		\hline
		$n=3$ &$2\cdot10^{-6}$& $1\cdot10^{-4}$ &$1$&\cellcolor{yellow} $3\cdot10^{7}$
		\\
		\hline
		$n=4$&$7\cdot10^{-5}$&  $1\cdot10^{-2}$ &$1\cdot10^{3}$&\cellcolor{yellow} $4\cdot10^{12}$\\
		\hline
		$n=5$ &$3\cdot10^{-3}$&$2$  & $8\cdot10^{5}$ &\cellcolor{yellow}$...$
		\\
		\hline
		$n=6$ &$0.1$& $246$  & \cellcolor{yellow}$4\cdot10^{8}$ &\cellcolor{yellow}$...$
		\\
		\hline
		$n=7$ &$13$& $2\cdot10^{4}$  & \cellcolor{yellow}$1\cdot10^{11}$ &\cellcolor{yellow}$...$
		\\
		\hline
		$n=8$ &$1\cdot10^{3}$&$5\cdot10^{6}$  & \cellcolor{yellow}$8\cdot10^{13}$ &\cellcolor{yellow}$...$
		\\
		\hline
		$n=9$ &$3\cdot10^{3}$&\cellcolor{yellow}$2\cdot10^{7}$  & \cellcolor{yellow}$9\cdot10^{14}$ &\cellcolor{yellow}$...$
		\\
		\hline
		$n=10$ &$1\cdot10^{5}$&\cellcolor{yellow}$1\cdot10^{9}$  & \cellcolor{yellow}$1\cdot10^{17}$ &\cellcolor{yellow}$...$
		\\
		\hline
		
	\end{tabular}\\
	\caption{The approximate time in seconds to compute all  $\hat{g}_i$ in the pumped-up solution of size $n^{2^k}$.   One year is $\approx 3.1 \cdot10^7$ seconds,  one month  is $\approx 2.6 \cdot10^6$ seconds and  one day is $\approx 8.6 \cdot10^4$ seconds.}\label{table_time_intruder}
\end{table}\\
In case, the original solution  has been previously pumped-up, we compute in Table \ref{table_time_intruder-more10},  the time required to   compute all  $\hat{g}_i$ in the pumped-up solutions of size $n^{2^k}$, where $n=m^2$, for $4\leq  m \leq 10$. In the computation, we consider the number of solutions for $m$, that is if $n=36$, then the time computed  is  $8.2 \cdot10^4$ times the time required to compute one $\hat{g}_i$ in this case. Clearly, this is far less time than the real time to compute $\hat{g}_i$  for all solutions of size $36$, but as we do not know all the solutions for $n>10$, this gives at least some estimation.

\begin{table}[h] 
	\centering
	\begin{tabular}{ |p{2cm}|p{2cm}|p{2cm}|p{2cm}| }
		\hline
		
		&$k=2$ &	$k=3$ &$ k=4$  \\
		\hline
		
		$n=16$ &$0.01$&$1\cdot10^{3}$  & \cellcolor{yellow}$4\cdot10^{12}$ 
		\\
		\hline
		$n=25$ &$2$&$8\cdot10^{5}$  &  \cellcolor{yellow}$...$ 
		\\
		\hline
		$n=36$ &$246$&\cellcolor{yellow}$4\cdot10^{8}$  & \cellcolor{yellow} $...$ 
		\\
		\hline
		$n=49$ &$2\cdot10^{4}$&\cellcolor{yellow}$1\cdot10^{11}$  & \cellcolor{yellow} $...$ 
		\\
		\hline
		$n=64$ &$5\cdot10^{6}$& \cellcolor{yellow}$8\cdot10^{13}$  & \cellcolor{yellow}$...$ 
		\\
		\hline
		$n=81$ &\cellcolor{yellow}$3 \cdot10^{7}$& \cellcolor{yellow}$1\cdot10^{15}$  & \cellcolor{yellow} $...$ 
		\\	\hline
		$n=100$ &\cellcolor{yellow}$1\cdot10^{9}$&  \cellcolor{yellow}$1\cdot10^{17}$  &  \cellcolor{yellow}$...$ 
		\\
		\hline
	\end{tabular}\\
	\caption{The approximate time in seconds to compute all  $\hat{g}_i$ in the pumped-up solution of size $n^{2^k}$, for an original solution which has been previously pumped-up.   One year is $\approx 3.1 \cdot10^7$ seconds,  one month  is $\approx 2.6 \cdot10^6$ seconds and  one day is $\approx 8.6 \cdot10^4$ seconds.}\label{table_time_intruder-more10}
\end{table}
From Tables \ref{table_time_one_gi}, \ref{table_time_intruder}, and \ref{table_time_intruder-more10}, there are several values of $n$ and $k$ for which the time required for Alice and Bob to compute the secret key is  short or reasonable while the time required for an intruder is  very long. As an example, the methods are secure for $n=9,10,36$ and $k=3$, or for $n=81,100$ and $k=2$ and some others. 
We need to recall that, as the intruder  does not know $n$, he needs to go through all  the possible values of $n$ until he breaks it.  So, the time required is the sum of  the times for all  these possible values of $n$,  and a given $k$.
\subsection{Some remarks on the strengths and weaknesses of the suggested methods }
\begin{enumerate}
	\item  In the encryption/decryption  method,  given the public key $i$, Alice knows the  encryption procedure of Bob (the permutation $\hat{g}_{i}\in S_{n^{2^k}}$) and Bob can  also compute the decryption procedure of Alice ($\hat{g}_{i}^{-1}\in S_{n^{2^k}}$). In the signature method, given $i$, Bob knows both the encryption and decryption procedures of Alice ($\hat{g}_i$  and $\hat{g}^{-1}_i$ respectively) and given $j$,  Alice  knows both the encryption and decryption procedures of Bob  ($\hat{g}_j$  and $\hat{g}^{-1}_j$). This is not the case in standard  procedures.
\item As Table \ref{table_time_one_gi} illustrates it, for $k \geq 5$, the time of computation of    $\hat{g}_i$ in the pumped-up solution of size $n^{2^k}$,  with the method suggested in Section  \ref{sec_tree}, is too large. We believe  there are much more efficient ways to do it  which may reduce the time of computation. Moreover,  for large $k$, choosing a small value for $i$  reduces also  the time of computation. The question is  how small this time of computation  can be. Indeed, we have to recall that already  $3^{2^5}\approx 1.8 \cdot10^{15}$, and $\hat{g}_i$ is  a permutation on such a very large number of elements. 

\item  As said above, the time of   computation of  the permutation $\hat{g}_i$ in the pumped-up solution of size $n^{2^k}$ might be large. Yet, it  is independent of the message, and so it can be done ahead of time if needed. As an example, for $n=49$ and $k=3$ (or $n=7$ and $k=4$), on one hand, the time required to compute $\hat{g}_i$  is relatively long, one day, but on the other hand it can be computed long before the sending of the  message. It might be worth considering  this kind of case, since from Table \ref{table_time_intruder-more10}, it may  take  thousands of  years for the intruder to break it.

\item Although the time of   computation of  the permutation $\hat{g}_i$ in the pumped-up solution of size $n^{2^k}$ might be large, it could be advisable to change  regularly  the public key $i$, and maybe also the other public key $k$.

	\item With the gap program, it is possible to have the list of all the solutions up to isomorphism, for $n\leq 10$. Clearly, it is more secure if Bob and Alice do not choose  their  original solution (their secret key) from the gap list, they can choose   an isomorphic solution to such a solution. Moreover, to enforce the security of the method, it could be advisable to use at each interaction another isomorphic solution to the original solution. Indeed, they can agree on some  permutation $\mu$ of large order  to apply on  the original solution to obtain  isomorphic solutions. From the proof of Lemma \ref{lem-equiv} $(ii)$,  there is no need to compute  $\hat{g}_i$ at each time, and it is enough to apply $\mu$ 
	on  $\hat{g}_i$ of the original solution (which as said above can be computed ahead of time).

	\item As $n$ grows, the total number of solutions of  size $n$ grows very fast and so for $n>10$, the number of solutions and their determination has not been established until now. We refer to \cite{arxiv-ybe} for details. The domain of possible values of $n$  is then restricted to the set $\{2,...,10, 16,25,36,49,64,81,100\}$. Clearly,   the security of the methods is enforced the  larger the domain of possible values of $n$ is.  So, it would be better to enlarge  the domain of $n$. From Tables  \ref{table_time_one_gi}-\ref{table_time_intruder-more10}, for $11 \leq n \leq 36$, and $k=3$, the time required to compute $\hat{g}_i$  is relatively short, while  it may  take  thousands of  years for the intruder to break it. So, for an implementation of this scheme, it   is not necessary to find all the solutions  for $11 \leq n \leq 36$, but it may be worth finding some  solutions of each size which have an interesting structure, in the sense that the cycle decompositions are diverse enough, and test their security.
	\item One may consider instead of a non-degenerate and involutive set-theoretic solution of the Yang-Baxter equation, a non-degenerate set-theoretic solution which is  non-involutive.  Such  a solution is  still defined by permutations. From the proof of Lemma \ref{lem-new--invol+braided}, the pumped-up solution is a non-degenerate set-theoretic solution which is  also non-involutive. The advantage  in considering non-involutive solutions is in their number. As an example, for $n=8$, there are   $34,530$  non-degenerate involutive solutions up to isomorphism, while there are $422,449,480$ non-degenerate non involutive solutions \cite{arxiv-ybe}. One needs to check that the suggested methods are still well defined in this case.
\end{enumerate}
To conclude, we have suggested a method for implementing a public-key 
cryptosystem and a method for signature, whose security, in both cases,  rests in  the difficulty  for an intruder to guess a permutation in a very large group of permutations.  Of course, the security of this method  needs to be examined in more detail. In particular, as the author is not a specialist in the domain of cryptography (or more precisely "a cryptographic amateur" as termed by W. Diffie and M.E.  Hellman  in \cite{diffie}), there may be some obvious weaknesses of the system that she is not aware of.
\section{Appendix: A tentative key exchange method}
 An ideal situation would be if there was a key exchange procedure in which the common secret key exchanged is the original solution. But, we do not know how to do that. Instead, we suggest a tentative key exchange method. It is tentative only, because it requires that both Alice and Bob share a common knowledge, the original solution, which is an inconvenient in such a procedure..\\
\textbf{Our suggestion of a key exchange  method:}\\

\begin{itemize}
	\item Begin with an original solution and pump it to get a solution of size $n^{2^k}$, for  $k$ large enough. 
	\item  Make $2^k$ (or $k$) public. Although  $2^k$  is  public, the size of the new solution, the size of the original solution   and the permutations defining the original solution are  effectively hidden from everyone else.
	\item   Choose $ 1 \leq i  \leq n^{2^k}$ such that  $i << n^{2^k}$ and make it public.
	\item  Bob chooses  $1 \leq j \leq n^{2^k}$, computes $\hat{g}_i(j)$ and sends it to Alice. 
	\item  Alice chooses  $1 \leq l \leq n^{2^k}$, computes $\hat{g}_i(l)$ and sends it to Bob.  
	\item Bob retrieves $l$ and computes the permutation  $\hat{g}_l\hat{g}_{g^{-1}_l(j)}$.
	\item  Alice retrieves $j$ and computes the permutation  $\hat{g}_j\hat{g}_{g^{-1}_j(l)}$.
	\item The common secret key is 	 $\hat{g}_l\hat{g}_{g^{-1}_l(j)}$	$=$ $\hat{g}_j\hat{g}_{g^{-1}_j(l)}$. Note that the  equality of these two permutations result  from Equation (\ref{eqn-braided-g--}) 
\end{itemize}
We illustrate with an  artificially  very simple example the key exchange method:  the original solution is from Example \ref{exemple-irretractable+indecomposable4} and it is pumped up    to a solution of size $n^{2^1}$, where $n=4$.  A part of the computations  is done in Example \ref{ex-pump4}. 
\begin{ex}
	The public key is  $2^1$ and make  $i=2$  public. Note that here $\hat{g}_i=\tilde{g}_i$. \\
	Bob chooses  $ j =3$, computes $\hat{g}_2(3)=6$ (since $g_1^2(T_1^3)=T_2^2$) and sends $6$ to Alice. \\
	Alice chooses  $l=14$, computes $\hat{g}_2(14)=1$ (since $g_1^2(T_4^2)=T_1^1$) and sends $1$  to Bob.  \\
	Bob retrieves $l=14$ and computes the permutation  $\hat{g}_{14}\hat{g}_{g^{-1}_{14}(3)}=\,\hat{g}_{14}\hat{g}_{4}$.\\
	Alice retrieves $j=3$ and computes the permutation  $\hat{g}_3\hat{g}_{g^{-1}_3(14)}=\,\hat{g}_3\hat{g}_{10}$.\\
	The common secret key is 	 $\hat{g}_3\hat{g}_{10}$	$=$ $\hat{g}_{14}\hat{g}_{4}$.  From Remark \ref{rem-renumber}, using the renumbering method,  this is equivalent to 
	$g_4^2g_1^4=g_1^3g_3^2$.
\end{ex}

\bigskip\bigskip\noindent
{ Fabienne Chouraqui}

\smallskip\noindent
University of Haifa at Oranim, Israel.

\smallskip\noindent
E-mail: {\tt fabienne.chouraqui@gmail.com} \\

                {\tt fchoura@sci.haifa.ac.il}
\end{document}